\documentclass[12pt]{article}           
\usepackage{amsmath}
\usepackage{amssymb}
\usepackage[mathscr]{eucal}
\usepackage{epsfig}
\usepackage{eufrak}
\newcommand{\por}[1]{\mbox{\boldmath${#1}$}}

\def\plb#1 #2 {Phys. Lett. {\bf #1B} #2 }
\def\phr#1 #2 {Phys. Rep. {\bf  #1} #2 }        
\def\npb#1 #2 {Nucl. Phys. {\bf B#1} #2 }
\def\aph#1 #2 {Ann. Phys. {\bf #1} #2 }         
\def\jmp#1 #2 {J. Math. Phys. {\bf #1} #2 }
\def\jgp#1 #2 {J. Geom. Phys. {\bf #1} #2 }
\def\prd#1 #2 {Phys. Rev. {\bf D#1} #2 }
\def\prl#1 #2 {Phys. Rev. Lett. {\bf #1} #2 }
\def\rmp#1 #2 {Rev. Mod. Phys.  {\bf #1} #2 }
\def\zpc#1 {Z. Phys. {\bf #1C} }
\def\cmp#1 #2 {Commun. Math. Phys. {\bf #1} #2 }
\def\cqg#1 #2 {Class.Quant.Grav. {\bf #1} #2 }
\def\mpl#1 {Mod. Phys. Lett. {\bf A#1} }
\def\cpc#1 {Computer Phys. Commun. {\bf #1} }   
\def\ijmp#1 {Int. J. Mod. Phys. {\bf A#1} }
\def\ijmpC#1 {Int. J. Mod. Phys. {\bf C#1} }

\begin{document}



\hfill
TUW--98--04

\vspace{40mm}

\begin{center}
{\huge\bf On the Cohomology of  Twisting  \\

\vspace{3mm}

Sheaves  on  Toric Varieties}
\end{center}

\vspace{10mm}

\begin{center} \vskip 12mm
  
{\large Mahmoud Nikbakht--Tehrani
\footnote[2]{e-mail: \texttt{nikbakht@tph16.tuwien.ac.at}}}

\vskip 5mm
       Institut f\"{u}r Theoretische Physik, Technische Universit\"{a}t Wien\\
       Wiedner Hauptstra\ss{}e 8--10, A-1040 Wien, AUSTRIA
\end{center}

\vspace{20mm}

\begin{center}
\begin{minipage}{12.6cm}
\begin{center}
{\bf Abstract}
\end{center}
 Using the homogeneous coordinate ring construction of a toric 
 variety ${\Bbb P}_{\scriptscriptstyle\Sigma}$  defined by a
 complete simplicial fan $\Sigma$ and the methods of local cohomology 
 theory we develop a framework for the calculation of cohomology 
 groups $H^{\bullet}({\Bbb P}_{\scriptscriptstyle\Sigma}, 
 {\cal O}_{{\Bbb P}_{\scriptscriptstyle\Sigma}}({\por p}))$ of 
 twisting sheaves ${\cal O}_{{\Bbb P}_{\scriptscriptstyle\Sigma}}
 ({\por p})$ on ${\Bbb P}_{\scriptscriptstyle\Sigma}$.  
\end{minipage}
\end{center}

\vspace{25mm}

February 1998
\thispagestyle{empty}


\newpage


\section{Introduction}\label{introduction}
\hrule 

\vspace{10mm}

The study of $(0,2)$ heterotic string compactifications leads us 
to consider holomorphic stable vector bundles \cite{wi86} (or more 
generally stable torsion-free or reflexive sheaves \cite{di96})
$\por{\mathscr E}$ on  Calabi-Yau varieties $X$. 
The physical properties of the resulting 
string vacua such as their stability and  massless particle 
spectrum are closely related to the topological properties 
of these bundles or sheaves.\\

In a large class of $(0,2)$ heterotic compactifications 
which has been constructed in the framework of gauged linear 
sigma models \cite{wi93,di94} the Calabi-Yau varieties $X$ are
realized as complete intersections of hypersurfaces in toric
varieties ${\Bbb P}_{\scriptscriptstyle\Sigma}$ and the  
relevant sheaves $\por{\mathscr E}$ are defined by short exact 
sequences or by the cohomology of monads in terms of the restriction 
of twisting sheaves ${\cal O}_{{\Bbb P}_{\scriptscriptstyle\Sigma}}
({\por p})$ to $X$. The knowledge of the cohomology 
groups of such sheaves in these models enables us to get  
useful information about the various topological invariants 
associated to $\por{\mathscr E}$ such as its cohomology groups 
and the like. The development of an efficient method for the 
calculation of such cohomology groups is the aim of the present work.\\

The basic ideas used here are those of the 
homogeneous coordinate ring construction of toric varieties 
defined by complete  simplicial fans \cite{cox92} and  local
cohomology theory \cite{gro67}. Combining these ideas with the
methods of computational commutative algebra we present an 
algorithmic way for approaching such cohomology calculations.\\

The outline of the paper is as follows. The next section is devoted
to a short survey of relevant concepts from the homogeneous coordinate
ring approach as they apply to our work. In section 3 we first review
the basic definitions and results on the sheaf cohomology and then 
go on to discuss its relation to  local cohomology theory. 
Afterwards we explain how these ideas apply to our situation.
In the last part of this section we describe some algorithms 
from commutative algebra
which will be used in the calculation of the local cohomology. 
In section 4, performing a few examples, we demonstrate explicitely 
the application of the developed methods. 
An appendix on the definition of the
inductive limit and of the Ext functor concludes this work.\\

\vspace{5mm}      
            

\section{Homogeneous coordinate ring approach}\label{homogeneous}
\hrule 

\vspace{10mm}

In this section we briefly discuss the relevant concepts from the 
homogeneous coordinate ring approach as they apply to our work.
The original motivation for its development  
was the desire to have a construction of toric varieties 
and related objects similar to those of ${\Bbb P}^{n}$ in classical 
algebraic geometry. For details on the constructions that follow  and 
for proofs we refer to \cite{cox92,DA78}.\\

To begin with we first introduce some notation. Let $\bf N$ and 
${\bf M}=\mbox{Hom}({\bf N},{\Bbb Z})$ denote a dual pair of lattices of
rank $d$ and $\langle \cdot , \cdot \rangle$ be the canonical pairing on
${\bf M}\times {\bf N}$. Further, let ${\bf N}_{\Bbb R}={\bf N}
\otimes_{\Bbb Z}{\Bbb R}$ and  ${\bf M}_{\Bbb R}={\bf M}
\otimes_{\Bbb Z}{\Bbb R}$ be the $\Bbb R$-scalar extensions of
$\bf N$ and $\bf M$, respectively. $T={\bf N}\otimes_{\Bbb Z} 
{\Bbb C}^{*}=\mbox{Hom}_{\Bbb Z}({\bf M},{\Bbb C}^{*})$ is
the $d$ dimensional algebraic torus which acts on the 
toric variety ${\Bbb P}_{\scriptscriptstyle \Sigma}$ defined by the  
(complete simplicial) fan $\Sigma$ in ${\bf N}_{\Bbb R}$ .
For a cone $\sigma \in \Sigma$ the dual cone, $\sigma^{\vee}$, is 
defined as usual by $\sigma^{\vee}=\{ m\in{\bf M}_{\Bbb R}\;|\; 
\langle m,n\rangle\geq 0 \;\; \mbox{for all} \;\; n\in \sigma \}$ 
and $\mbox{cosp}\sigma^{\vee}$ is the greatest subspace of 
${\bf M}_{\Bbb R}$ contained in $\sigma^{\vee}$. The open affine 
variety in ${\Bbb P}_{\scriptscriptstyle \Sigma}$ associated to $\sigma$ 
is denoted by $X_{\sigma^{\vee}}$. Let $\Sigma^{(k)}$ be the set 
of $k$ dimensional cones in $\Sigma$ . By $e_{i}$ we denote the 
primitive lattice vectors on the one dimensional cones in 
$\Sigma^{(1)}=\{\rho_{1}, \ldots,\rho_{n}\}$. This set will play an 
important role in what follows.\\

Each one dimensional cone $\rho_{i}$ defines a $T$-invariant Weil divisor,
denoted by $D_{i}$, which is the closed subvariety 
$X_{\mbox{\scriptsize cosp}\rho_{i}^{\vee}}$ in $X_{\rho_{i}^{\vee}}$. 
This is indeed the closed $T$-orbit associated to $\rho_{i}\,$. 
The finitely generated free abelian group $\bigoplus_{i=1}^{n}{\Bbb Z}
\cdot D_{i}$ is the group of $T$-invariant Weil divisors in 
${\Bbb P}_{\scriptscriptstyle \Sigma}$. Each $m\in {\bf M}$ gives a character 
$\chi^{m}:T\to {\Bbb C}^{*}$, and hence $\chi^{m}$ is a rational 
function on ${\Bbb P}_{\scriptscriptstyle \Sigma}$. It defines
the $T$-invariant Cartier divisor 
$\mbox{div}(\chi^{m})=\sum_{i=1}^{n}\langle m,e_{i}
\rangle \; D_{i}\;$~. In this way we obtain the map {\boldmath$\alpha$}
\begin{eqnarray}
\mbox{\boldmath$\alpha$}: {\bf M}\longrightarrow 
\bigoplus_{i=1}^{n}{\Bbb Z}\cdot D_{i}\;\; ,
\;\; m\mapsto \sum_{i=1}^{n}\langle m,e_{i}\rangle \; D_{i}\;\label{2.1} .
\end{eqnarray}
It follows from the completeness of the fan
$\Sigma$ that the map {\boldmath$\alpha$} is injective.
The cokernel of this map defines the Chow group $A_{d-1}
({\Bbb P}_{\scriptscriptstyle \Sigma})$ which is a finitely generated abelian group
of rank $n-d\,$. Therefore we have the following
exact sequence
\begin{eqnarray}
0\longrightarrow {\bf M}\stackrel{\mbox{\boldmath$\alpha$}}{\longrightarrow}
\bigoplus_{i=1}^{n}{\Bbb Z}\cdot D_{i}\stackrel{deg}{\longrightarrow}
A_{d-1}({\Bbb P}_{\scriptscriptstyle \Sigma})\longrightarrow 0\;\label{2.2},
\end{eqnarray}
where $\mbox{\sl deg}$ denotes the canonical projection. Now consider 
$G=\mbox{Hom}_{\Bbb Z}(A_{d-1}({\Bbb P}_{\scriptscriptstyle \Sigma}),{\Bbb C}^{*})$
which is in general isomorphic to a product of 
$({\Bbb C}^{*})^{n-d}$ and a finite group.
By applying $\mbox{Hom}_{\Bbb Z}(\;\cdot \;, {\Bbb C}^{*})$ to 
(\ref{2.2}) we get
\begin{eqnarray}
1\longrightarrow G \longrightarrow ({\Bbb C}^{*})^{n}
\longrightarrow T\longrightarrow 1\; \nonumber,
\end{eqnarray}
which defines the action of $G$ on ${\Bbb C}^{n}$ : 
\begin{eqnarray}
g\cdot (x_{1},\ldots ,x_{n})=(\,g( \mbox{\sl deg}\,D_{1})\; x_{1},\ldots ,
                               g( \mbox{\sl deg}\,D_{n})\; x_{n}\,) \nonumber 
\end{eqnarray}
for $g\in G$ and $(x_{1},\ldots,x_{n})\in {\Bbb C}^{n}$.\\

Let $S={\Bbb C}[x_{1},\ldots,x_{n}]$ be the polynomial ring over $\Bbb C$
in the variables $x_{1},\ldots,x_{n}$, where  the $x_{i}$ correspond to the
one-dimensional cones $\rho_{i}$ in $\Sigma$ . 
Each monomial $\,x_{1}^{a_{1}}\ldots x_{n}^{a_{n}}\,$  in $S$ determines a
divisor $\,\sum_{i=1}^{n}a_{i}\,D_{i}\,$. This ring is graded 
in a natural way by $\,\mbox{\sl deg}(x_{i}):=\mbox{\sl deg}\,D_{i}\,$:
\begin{eqnarray}
    S=\hspace{-3mm}\bigoplus_{{\por {\scriptscriptstyle q}}
      \in A_{d-1}({\Bbb P}_{\scriptscriptstyle\Sigma})}
      \hspace{-2mm}S_{\por {\scriptscriptstyle q}}\, ,\nonumber
\end{eqnarray}
where $S_{\por{\scriptscriptstyle q}}$ is generated by all  monomials $\,x_{1}^{a_{1}}\ldots 
x_{n}^{a_{n}}\,$  such that $\mbox{\sl deg}\,(\sum_{i=1}^{n}
a_{i}D_{i})={\por q}\,$. Let $I$ denote the monomial ideal in 
$S$ generated by ${\por x}^{\sigma}=\prod_{\rho_{i}\not\prec \sigma }
\,x_{i}$ for all $\sigma\in\Sigma\,$. Note that the set of monomials 
$\{ {\por x}^{\sigma}\mid \,\sigma\in\Sigma^{(d)}\,\}\,$, 
is the unique minimal basis of this ideal. 
The ring $S$ defines the $n$-dimensional 
affine space ${\Bbb A}^{n}= \mbox{Spec}\,S$. The ideal $I$ gives 
the variety 
$\;{\por Z}_{\scriptscriptstyle\Sigma}={\por V}(I)\;$
which is denoted as the exceptional set. 
Removing the exceptional set ${\por Z}_{\scriptscriptstyle\Sigma}$ we obtain
the Zariski open set
$\;{\por U}_{\scriptscriptstyle \Sigma}= {\Bbb A}^{n}\setminus
                       {\por Z}_{\scriptscriptstyle \Sigma}\;$ 
which is invariant under the action of $G$.
For the case of a complete simplicial fan  the geometric
quotient of ${\por U}_{\scriptscriptstyle \Sigma}$ by $G$ 
exists and gives rise to 
${\Bbb P}_{\scriptscriptstyle \Sigma}$ \cite{cox92}.\\

Using this construction of a toric variety 
${\Bbb P}_{\scriptscriptstyle\Sigma}$,
the sheaves of ${\cal O}_{{\Bbb P}_{\scriptscriptstyle\Sigma}}$-modules can be
studied in a way similar to that of ${\Bbb P}^{n}$. For example,
the twisting sheaves ${\cal O}_{{\Bbb P}_{\scriptscriptstyle\Sigma}}({\por p})$
are sheaves which are associated to the graded $S$-module $S({\por p})$,
i.e. ${\cal O}_{{\Bbb P}_{\scriptscriptstyle\Sigma}}({\por p})= S({\por p})^{\sim}$,
where $S({\por p})_{\por {\scriptscriptstyle q}}=S_{{\por {\scriptscriptstyle p}}+{\por {\scriptscriptstyle q}}}$
for all ${\por q}\,$.  
It can be shown that 
${\cal O}_{{\Bbb P}_{\scriptscriptstyle\Sigma}}({\por p})\simeq  
{\cal O}_{{\Bbb P}_{\scriptscriptstyle\Sigma}}(D)$, 
where $D$ is a $T$-invariant
divisor with $\mbox{\sl deg}\,D={\por p}$. In particular, 
${\cal O}_{{\Bbb P}_{\scriptscriptstyle\Sigma}}=S^{\sim}$. Furthermore, we
have $S\simeq \bigoplus_{\por {\scriptscriptstyle p}}H^{0}({\Bbb P}_{\scriptscriptstyle\Sigma},
{\cal O}_{{\Bbb P}_{\scriptscriptstyle\Sigma}}({\por p}))$. We conclude 
this section with the following theorem \cite{cox92}.\\

\noindent
{\bf Theorem A}: Let $\Sigma$ be a complete simplicial fan. Then every
coherent sheaf ${\cal F}$ on ${\Bbb P}_{\scriptscriptstyle\Sigma}$ is of the form
${\cal F}=M^{\sim}$, where $M$ is a finitely generated graded
$S$-module.\\



\section{Cohomology of twisting sheaves}\label{cohomology}
\hrule

\vspace{10mm}

We begin this section with recalling some definitions and facts from 
the cohomology of sheaves. For more details on the following 
concepts and the proofs we refer to the
standard works \cite{gro67,ser55,cox921}.\\

Let ${\cal F}$ be an Abelian sheaf on 
a topological space $X$. Consider the sheaf $C^{0}({\cal F})$ defined by
\begin{eqnarray}
     C^{0}({\cal F})\,(U)=\prod_{x\in U}{\cal F}_{x}\;.\nonumber    
\end{eqnarray}
The sheaf ${\cal F}$ is canonically embedded in $C^{0}({\cal F})$
by associating to $s\in {\cal F}(U)$ the family $(s(x))\in 
\prod_{x\in U}{\cal F}_{x}\,$. The sheaf  $C^{0}({\cal F})$
is always flabby. Now consider the family $\{C^{n}({\cal F})\}_{n\geq 1}$
of sheaves that is recursively defined by
\begin{eqnarray}
 & & C^{1}({\cal F}):=C^{0}(C^{0}({\cal F})\,/\,{\cal F})\, , \nonumber\\
 & & C^{n+1}({\cal F}):=C^{0}(C^{n}({\cal F})\,/\,
                        d^{n-1}\,C^{n-1}({\cal F}))\, ,\nonumber
\end{eqnarray}
where $d^{n}:C^{n}({\cal F})\to C^{n+1}({\cal F})$ is defined as a composition
\begin{eqnarray}
   C^{n}({\cal F})  \longrightarrow 
   (\,C^{n}({\cal F})\,/\, d^{n-1}\,C^{n-1}({\cal F})\,)\longrightarrow
   C^{0}(\,C^{n}({\cal F})\,/\, d^{n-1}\,C^{n-1}({\cal F})\,)\; .\nonumber
\end{eqnarray}
In this way we obtain the Godement canonical flabby resolution
$0\to {\cal F}\to C^{\bullet}({\cal F})\,$. The cohomology of the
complex $C^{\bullet}({\cal F})\,(X)$ of Abelian groups is said to be
the cohomology of the sheaf ${\cal F}$ and is 
denoted by $H^{\bullet}(X,{\cal F})$.
Clearly $\Gamma (X,{\cal F})=H^{0}(X,{\cal F})$ and 
$H^{p}(X,{\cal F})=0$ for $p<0$. If ${\cal F}$ is flabby, then
$H^{p}(X,{\cal F})=0$ for $p>0$. The definition of the cohomology
group given here is not very enlightening. For this reason 
we now discuss some other methods which provide us with more useful
tools for the cohomology computations.\\

First, we recall the definition of the \v{C}ech cohomology group. Let
${\mathscr U}=\{ U_{i} \}_{i\in I}$ be an open covering of $X$,
where the index set is well-ordered. For an Abelian sheaf ${\cal F}$
on $X$ we set 
\begin{eqnarray} 
   C^{p}({\mathscr U},{\cal F}):=\bigoplus_{i_{0}< \ldots < i_{p}}
   \Gamma (U_{i_{0}}\cap\ldots\cap  U_{i_{p}},{\cal F}) \nonumber
\end{eqnarray}
and define $d^{p}:C^{p}({\mathscr U},{\cal F})\to C^{p+1}
({\mathscr U},{\cal F})$ by
\begin{eqnarray} 
   (d^{p}\,\alpha)_{i_{0}\ldots i_{p+1}}=\sum_{k=0}^{p+1}
   (-1)^{k}\,{\alpha_{i_{0}\ldots \hat{i}_{k}\ldots i_{p+1}}}_
   {\mid_{U_{i_{0}}\cap \ldots \cap U_{i_{p+1}}}}\; ,\nonumber
\end{eqnarray}
where the $\;\hat{}\;$ over $i_{k}$ means that this index is to be omitted.
$C^{\bullet}({\mathscr U},{\cal F})$ is called the \v{C}ech complex
for the open covering ${\mathscr U}$ with values in ${\cal F}$.
The cohomology groups $\check{H}^{\bullet}({\mathscr U},{\cal F})$ of
this complex are called the \v{C}ech cohomology groups of ${\cal F}$
for the open covering ${\mathscr U}$.\\

We are primarily concerned with the coherent sheaves 
(of ${\cal O}_{X}$-modules) on an algebraic variety
$(X,{\cal O}_{X})$. The following theorem of Serre plays a fundamental
role in the cohomology theory of such sheaves.\\

\noindent
{\bf Theorem (Serre)}: Let ${\cal F}$ be a quasi-coherent sheaf on an
affine variety $X$. Then $H^{p}(X,{\cal F})=0$ for $p>0\,$.\\

If the variety $X$ is separated, then the following theorem 
gives us a useful tool for the computation of cohomology groups.\\

\noindent
{\bf Theorem B}: Let ${\cal F}$ be a quasi-coherent 
sheaf on a separated  variety $X$ and ${\mathscr U}=\{U_{i}\}_{i=0}^{n}$
be an open covering of $X$ by open affine subvarieties. Then
$\check{H}^{p}({\mathscr U},{\cal F})\simeq H^{p}(X,{\cal F})$.\\

The key point in the proof of this theorem is to
show that $H^{q}(U_{i_{0}}\cap\ldots\cap U_{i_{p}},{\cal F})=0$
for all $(p+1)$-tuples $i_{0}<\ldots <i_{p}$ and $q>0$. 
This follows already from Serre's theorem because the separability of 
$X$ implies that $ U_{i_{0}}\cap\ldots\cap U_{i_{p}}$ are affine.\\

Having reduced the calculation of cohomology to
that of \v{C}ech cohomology for a finite covering, we now come to 
the question how these cohomology groups can be effectively
calculated. The Koszul complex method provides us with the 
necessary tools. Before going further, we mention the following
vanishing theorem of Grothendieck.\\

\noindent
{\bf Theorem (Grothendieck)}: Let $X$ be an $n$-dimensional  Noetherian 
topological space. Then, $H^{p}(X,{\cal F})=0$ for all $p>n$ and 
all Abelian sheaves ${\cal F}$ on $X$.\\  

Now, let $R$ be a commutative ring with unit. Consider the free $R$-module
$R^{d}$ with the canonical basis $\{e_{i}\}_{i=1}^{d}\,$. Further, let
   $ d_{1}:R^{d}\longrightarrow R\, ,\,
    e_{i}\mapsto a_{i}:=d_{1}(e_{i})$
be a homomorphism of $R$-modules. We are going to construct a complex
of $R$-modules. Let 
\begin{eqnarray}
 K_{p}({\por a})=K_{p}(a_{1},\ldots,a_{d})=
 \bigwedge\nolimits^{p}R^{d}=\underset{i_{1}<\ldots < i_{p}}{\bigoplus}\,
 R\, e_{i_{1}}\wedge\ldots\wedge  e_{i_{p}}\,.\nonumber 
\end{eqnarray}
Define
$d_{p}:K_{p}({\por a})\to K_{p-1}({\por a})$ by
\begin{eqnarray}
   d_{p}\,(e_{i_{1}}\wedge\ldots\wedge  e_{i_{p}})=
   \sum_{k=1}^{p}(-1)^{k}d_{1}(e_{i_{k}})\,
   e_{i_{1}}\wedge\ldots\wedge \hat{e}_{i_{k}}
   \wedge\ldots\wedge  e_{i_{p}}\;.\nonumber
\end{eqnarray}
It can be easily seen that $d_{p-1}\circ d_{p}=0\,$. Therefore, we get the
following complex
\begin{eqnarray}
 0\to\bigwedge\nolimits^{d}R^{d}\to\bigwedge\nolimits^{d-1}R^{d}
 \to\ldots\to R^{d}\stackrel{d_{1}}{\longrightarrow}R\to 0\, ,\nonumber
\end{eqnarray}
which is called a (homological) Koszul complex. Let $M$ be a finitely
generated  $R$-module.
Then $K_{\bullet}({\por a},M):=K_{\bullet}({\por a})\otimes_{\scriptscriptstyle R}M$
(resp. $K^{\bullet}({\por a},M):=\mbox{Hom}_{\scriptscriptstyle R}(K_{\bullet}({\por a}),M)$)
is called the homological (resp. cohomological) Koszul complex of $M$ with
respect to $(a_{1},\ldots,a_{d})$. Its corresponding homology
(cohomology) groups are denoted by $H_{\bullet}({\por a},M)$
($H^{\bullet}({\por a},M)$). Note that $K^{\bullet}({\por a},M)
=\underset{i_{1}<\ldots < i_{p}}{\bigoplus}\,M\;$. Therefore, $\alpha
\in K^{\bullet}({\por a},M)$ can be considered as a sequence
$\alpha =(m_{i_{1}\ldots i_{p}})_{\;i_{1}<\ldots < i_{p}}$ of elements
of $M$. The coboundary operator $d^{p}:K^{p}({\por a},M)\to 
K^{p+1}({\por a},M)$ is given by
\begin{eqnarray} 
  (d^{p}\,\alpha)_{\,i_{1}\ldots i_{p+1}}=\sum_{k=1}^{p+1}(-1)^{k}\,
  d_{1}(e_{i_{k}})\;m_{i_{1}\ldots \hat{i}_{k}\ldots i_{p+1}}\; .\nonumber
\end{eqnarray}
It is not hard to see that $H_{0}({\por a},M)=M\,/\,
\langle a_{1},\ldots,a_{d}\rangle\,M$ and $H_{d}({\por a},M)
=\{\,m\in M\,\mid \, \langle a_{1},\ldots,a_{d}\rangle\, m =0\,\}\,$.\\

\noindent
{\bf Theorem C}: If $\,\langle a_{1},\ldots,a_{d}\rangle =R\,$, then 
$\,K_{\bullet}({\por a},M)\,$ is acyclic.\\

A sequence $\,(b_{1},\ldots,b_{r})\,$ in $R\,$ is said to be $M$-regular
if $\,\langle b_{1},\ldots,b_{r}\rangle \neq R\,$ and the image 
of $\,b_{i}\,$ in $\,M\,/\,\langle b_{1},\ldots,b_{i-1}\rangle\,M\,$ is
no zerodivisor for $\,i=1,\ldots,r\,$.\\

\noindent
{\bf Theorem D}: If $H_{p}({\por a},M)=0$ for $p>r$ while 
$H_{r}({\por a},M)\neq0$, then every maximal $M$-sequence in
$\langle a_{1},\ldots,a_{d}\rangle$ has length $r$. In particular,
if $\,(a_{1},\ldots,a_{d})\,$ is an $M$-regular sequence
in $R\,$, then $\,H_{0}({\por a},M)=M\,/\,\langle a_{1},
\ldots,a_{d}\rangle\,M\,$ and $\,H_{p}({\por a},M)=0\,$ for $\,p=1,
\ldots, d\,$.\\

\noindent
{\bf Theorem E}: $H_{p}({\por a},M)\simeq H^{d-p}({\por a},M)$.\\

After these preliminaries, we now turn our attention to the relation
between Koszul and \v{C}ech complexes. We will show that the 
\v{C}ech complex for a given covering is related to the limit
of an inductive system of Koszul complexes. Let $U_{i}$ be
the principal open subsets in $\mbox{Spec}\,R$ associated to
$a_{i}$ and ${\mathscr U}= \{U_{i}\}_{i=1}^{d}\,$. Let 
$U=\bigcup_{i=1}^{d}\,U_{i}$ and  ${\cal F}
=M^{\sim}\,$. Then $\Gamma (U_{i_{1}}\cap\ldots\cap U_{i_{p}},{\cal F})$
is equal to the localization
of $M$ at $a_{i_{1}}\ldots a_{i_{p}}$, i.e. $M_{\,a_{i_{1}}\ldots 
a_{i_{p}}}\,$. Now, consider the family $\{K^{\bullet}
({\por a}^{m},M)\}_{m\in {\Bbb N}}$ of Koszul complexes of $M$
with respect to ${\por a}^{m}:=(a_{1}^{m},\ldots,a_{d}^{m})\,$.
Note that $d^{p}_{m}: K^{p}({\por a}^{m},M)\to K^{p+1}({\por a}^{m},M)$
is given by 
\begin{eqnarray}
(d^{p}_{m}\,\alpha)_{i_{1}\ldots i_{p+1}}=\sum_{k=1}^{p+1}
\;(-1)^{k}\,a_{i_{k}}^{m}\;{m_{\,i_{1}\ldots \hat{i}_{k}
\ldots i_{p+1}}}\,.\nonumber 
\end{eqnarray}
Let $\,f^{p}_{mn}:K^{p}({\por a}^{m},M)\to
K^{p}({\por a}^{n},M)$ be defined by
\begin{eqnarray}
   f^{p}_{mn}:(m_{i_{1}\ldots i_{p}})_{\;i_{1}<\ldots < i_{p}}\mapsto
              (a_{i_{1}}\ldots a_{i_{p}})^{n-m}\; (m_{i_{1}\ldots 
               i_{p}})_{\;i_{1}<\ldots < i_{p}}\nonumber           
\end{eqnarray}
for $n\geq m$. It can be easily checked that $f^{p}_{mm}= \mbox{id}$ and
$f^{p}_{mn}\circ f^{p}_{\ell m}=f^{p}_{\ell n}\,$. Using the definitions
of coboundary operators and $f^{p}_{mn}$ one can show that
the following diagram commutes\\
\begin{center}
\unitlength0.25cm
\begin{picture}(48,14)
   \put(13.7,12){\mbox{$K^{p}({\por a}^{m},M)$}}
   \put(22,12.5){\vector(1,0){6}}
   \put(28.5,12){\mbox{$K^{p+1}({\por a}^{m},M)$}}
   \put(14,2){\mbox{$K^{p}({\por a}^{n},M)$}}
   \put(22,2.5){\vector(1,0){6}}
   \put(28.3,2){\mbox{$K^{p+1}({\por a}^{n},M)$}}
   \put(17,11){\vector(0,-1){7}}
   \put(33,11){\vector(0,-1){7}}
\end{picture}
\end{center}

Therefore, the family $\{K^{\bullet}({\por a}^{m},M)\}_{m\in {\Bbb N}}$
together with $\{f^{\bullet}_{mn}\}_{\,m\leq n}$ build an inductive
system of Koszul complexes. Let $K^{\bullet}({\por a}^{\infty},M)=
\underset{\longrightarrow}{\mbox{lim}}\,K^{\bullet}({\por a}^{m},M)\,$.\\

We now define the map $\,\varphi^{p}_{m}:K^{p}({\por a}^{m},M)\to
C^{p-1}({\mathscr U},{\cal F})\,$ by
\begin{eqnarray}
 (m_{\,i_{1}\ldots i_{p}})\mapsto 
 \left(\frac{(m_{\,i_{1}\ldots i_{p}})}{(a_{i_{1}}
        \ldots a_{i_{p}})^{m}}\right) \;.\nonumber
\end{eqnarray}
Obviously, the following diagram 
\begin{center}
\unitlength0.25cm
\begin{picture}(20,15)
   \put(1.3,12){\mbox{$K^{p}({\por a}^{m},M)$}}
   \put(12,11){\mbox{$\scriptstyle\varphi^{p}_{m}$}}
   \put(9.5,12.3){\vector(1,-1){4}}
   \put(14,7){\mbox{$C^{p-1}({\mathscr U},{\cal F})$}}
   \put(1.5,2){\mbox{$K^{p}({\por a}^{n},M)$}}
   \put(12,3){\mbox{$\scriptstyle\varphi^{p}_{n}$}}
   \put(9.5,2.7){\vector(1,1){4}}
   \put(5.5,11){\vector(0,-1){7}}
   \put(2.5,8){\mbox{$\scriptstyle f^{p}_{mn}$}}
\end{picture}
\end{center}
commutes for all $m\leq n$. Therefore, $C^{p-1}({\mathscr U},{\cal F})$ 
is a target object for the inductive system 
$\{K^{p}({\por a}^{m},M)\}_{m\in {\Bbb N}}\,$. Because of the universal
property of the inductive limit there exists a map $\psi^{p}:
K^{p}({\por a}^{\infty},M)\to C^{p-1}({\mathscr U},{\cal F})$ such
that the diagram
\begin{center}
\unitlength0.25cm
\begin{picture}(45,20)
   \put(2,17){\mbox{$K^{p}({\por a}^{m},M)$}}
   \put(13.4,14.6){\mbox{$\scriptstyle f^{p}_{m}$}}
   \put(10,17.3){\vector(1,-1){6}}
   \put(26,14.3){\mbox{$\scriptstyle\varphi^{p}_{m}$}}
   \put(10,17.7){\vector(4,-1){29.5}}
   \put(11.5,9.5){\mbox{$K^{p}({\por a}^{\infty},M)$}}
   \put(26,10.5){\mbox{$\scriptstyle\psi^{p}$}}
   \put(20,10){\vector(1,0){19.5}}
   \put(39.7,9.5){\mbox{$C^{p-1}({\mathscr U},{\cal F})$}}
   \put(2.2,2){\mbox{$K^{p}({\por a}^{n},M)$}}
   \put(13.4,5){\mbox{$\scriptstyle f^{p}_{n}$}}
   \put(10,2.7){\vector(1,1){6}}
   \put(26,5){\mbox{$\scriptstyle\varphi^{p}_{n}$}}
   \put(10,2.3){\vector(4,1){29.5}}
   \put(6,16){\vector(0,-1){12}}
   \put(3,10){\mbox{$\scriptstyle f^{p}_{mn}$}}
\end{picture}
\end{center}
commutes. We now prove that $\psi^{p}$ is actually an isomorphism.
Let $\beta\in C^{p-1}({\mathscr U},{\cal F})\,$. Then, for 
$i_{1}<\ldots <i_{p}\,$, $\beta_{i_{1}\ldots i_{p}}= (a_{i_{1}}\ldots 
a_{i_{p}})^{-m}\cdot m_{i_{1}\ldots i_{p}}$, where
$m_{i_{1}\ldots i_{p}}\in M\,$. Therefore, 
$\alpha=(m_{i_{1}\ldots i_{p}})_{\,i_{1}<\ldots <i_{p}}\in
K^{p}({\por a}^{m},M)$ and $\beta =\psi^{p}(f^{p}_{m}(\alpha ))\,$.
This shows the surjectivity of $\psi^{p}$. Injectivity of
$\psi^{p}$ : Let $f^{p}_{m}(\alpha )\in\mbox{Ker}\,\psi^{p}$, i.e.
$\psi^{p}(f^{p}_{m}(\alpha ))=0\,$. It follows from the commutativity
of the above diagram that $\psi^{p}(f^{p}_{m}(\alpha ))=(\psi^{p}\circ 
f^{p}_{m})\,(\alpha )=\varphi^{p}_{m}(\alpha)=0\,$, i.e.
$(a_{i_{1}}\ldots a_{i_{p}})^{-m}\cdot m_{i_{1}\ldots i_{p}}=0\,$ for
each $i_{1}<\ldots <i_{p}\,$. Consequently, there exists some $k$
such that $(a_{i_{1}}\ldots a_{i_{p}})^{k}\cdot m_{i_{1}\ldots i_{p}}=0\,$.
But $(a_{i_{1}}\ldots a_{i_{p}})^{k}\cdot m_{i_{1}\ldots i_{p}}=
(\,f^{p}_{m,m+k}(\alpha )\,)_{i_{1}\ldots i_{p}}\,$. That is 
$f^{p}_{m,m+k}(\alpha )=0\,$. On the other hand, $f^{p}_{m}(\alpha )=
(f^{p}_{m+k}\circ f^{p}_{m,m+k})\,(\alpha )=f^{p}_{m+k}(\,
f^{p}_{m,m+k}(\alpha )\,)=0\,$. Hence $\mbox{Ker}\,\psi^{p}=
\{0\}\,$.\\

Summarizing the above discussion, we can write down the following
exact sequence of complexes
\begin{center}
\unitlength0.25cm
\begin{picture}(48,20)
\put(-1.5,17){\mbox{$0$}}
   \put(16.5,17){\mbox{$0$}}
   \put(28.5,17){\mbox{$0$}}
   \put(39,17){\mbox{$0$}}
   \put(-4.3,12){\mbox{$\tilde{C}^{\bullet}({\mathscr U},{\cal F})$}}
   \put(5,12){\mbox{$:$}}
   \put(7,12){\mbox{$0$}}
   \put(8.5,12.5){\vector(1,0){6.7}}
   \put(16.5,12){\mbox{$0$}}
   \put(18.7,12.5){\vector(1,0){6}}
   \put(25.5,12){\mbox{$C^{0}({\mathscr U},{\cal F})$}}
   \put(32.5,12.5){\vector(1,0){4}}
   \put(37,12){\mbox{$C^{1}({\mathscr U},{\cal F})$}}
   \put(44.5,12.5){\vector(1,0){2}}
   \put(47,12.5){\mbox{$\ldots$}}
   \put(-5.2,7){\mbox{$K^{\bullet}({\por a}^{\infty},M)$}}
   \put(5,7){\mbox{$:$}}
   \put(7,7){\mbox{$0$}}
   \put(8.5,7.5){\vector(1,0){3}}   
   \put(12.5,7){\mbox{$K^{0}({\por a}^{\infty},M)$}}
   \put(21.7,7.5){\vector(1,0){2.5}}
   \put(25,7){\mbox{$K^{1}({\por a}^{\infty},M)$}}
   \put(33.7,7.5){\vector(1,0){2}}
   \put(36.5,7){\mbox{$K^{2}({\por a}^{\infty},M)$}}
   \put(45,7.5){\vector(1,0){1.7}}
   \put(47.2,7.5){\mbox{$\ldots$}}
   \put(-3,2){\mbox{$M[0]$}}
   \put(5,2){\mbox{$:$}}
   \put(7,2){\mbox{$0$}}
   \put(8.5,2.5){\vector(1,0){3}}   
   \put(12.5,2){\mbox{$K^{0}({\por a}^{\infty},M)$}}
   \put(21.7,2.5){\vector(1,0){6}}
   \put(28.5,2){\mbox{$0$}}
   \put(30.5,2.5){\vector(1,0){7}}
   \put(39,2){\mbox{$0$}}
   \put(41,2.5){\vector(1,0){5.8}}
   \put(47,2.5){\mbox{$\ldots$}}
\put(-1.5,-3){\mbox{$0$}}   
   \put(16.5,-3){\mbox{$0$}}
   \put(28.5,-3){\mbox{$0$}}
   \put(39,-3){\mbox{$0$}}
\put(-1.1,15.9){\vector(0,-1){2}}
   \put(17,15.9){\vector(0,-1){2}}
   \put(29,15.9){\vector(0,-1){2}}
   \put(39.5,15.9){\vector(0,-1){2}}   
\put(-1.1,10.8){\vector(0,-1){2}}
   \put(17,10.8){\vector(0,-1){2}}
   \put(29,10.8){\vector(0,-1){2}}
   \put(39.5,10.8){\vector(0,-1){2}}
\put(-1.1,5.8){\vector(0,-1){2}}
   \put(17,5.8){\vector(0,-1){2}}
   \put(29,5.8){\vector(0,-1){2}}
   \put(39.5,5.8){\vector(0,-1){2}}
\put(-1.1,0.8){\vector(0,-1){2}}
   \put(17,0.8){\vector(0,-1){2}}
   \put(29,0.8){\vector(0,-1){2}}
   \put(39.5,0.8){\vector(0,-1){2}}
\end{picture}
\end{center}\vspace{10mm}
which yields
\begin{eqnarray}
& & 0\to H^{0}({\por a}^{\infty},M)\to M\to \check{H}^{0}({\mathscr U},
    {\cal F}_{\mid U})\to H^{1}({\por a}^{\infty},M)\to 0 \nonumber\\
& & \label{ck}\\
& & \check{H}^{p-1}({\mathscr U},{\cal F}_{\mid U})\simeq 
    H^{p}({\por a}^{\infty},M)\hspace{8mm}\mbox{for}\;\; p>1\;.\nonumber
\end{eqnarray}

It is noteworthy that there is a close relation between
$\,H^{p}({\por a}^{\infty},M)\,$ and $\,\mbox{Ext}^{p}_{R}\,$ 
providing us with yet another useful computatianl tool: 
\begin{eqnarray}
   H^{p}({\por a}^{\infty},M)\simeq\underset{\longrightarrow}{\mbox{lim}}
   \;\mbox{Ext}^{p}_{R}(R\,/\,I^{m} , M)\; , \nonumber
\end{eqnarray}
where $I^{m}=\langle a_{1}^{m},\dots, a_{d}^{m} \rangle\,$. For the
details of the proof we refer to \cite{gro67} \S 2.\\

Now, the calculation of the cohomology of twisting sheaves on 
${\Bbb P}_{\scriptscriptstyle\Sigma}$ can be approached in the following way.
As mentioned before, in the homogeneous coordinate ring construction 
the toric variety ${\Bbb P}_{\scriptscriptstyle\Sigma}$ is realized as the geometric 
quotient ${\por U}_{\scriptscriptstyle\Sigma}\,/\, G\,$ provided that $\Sigma$ 
is a complete simplicial fan. An open covering of 
${\por U}_{\scriptscriptstyle\Sigma}$ is given by the
principal open subsets $U_{i}:=D({\por x}^{\sigma_{i}})$, where 
the ${\por x}^{\sigma_{i}}$ belong to the minimal basis of the ideal
$I\,$:
\begin{eqnarray} 
   {\por U}_{\scriptscriptstyle\Sigma}=\bigcup_{\sigma_{i}\in\Sigma^{(d)}} 
                        U_{i}\;.\nonumber
\end{eqnarray}
Note that the projection map $\,\pi : {\por U}_{\scriptscriptstyle\Sigma}\to
{\Bbb P}_{\scriptscriptstyle\Sigma}\,$ is an affine morphism. For such morphisms
we have the following lemma.\\

\noindent
{\bf Lemma}: Let $f:X\to Y$ be an affine morphism of separated varieties
$X$ and $Y$. Let ${\cal F}$ be a quasi-coherent ${\cal O}_{X}$-module. 
Then $H^{p}(X,{\cal F})\simeq H^{p}(Y,f_{*}{\cal F})$.\\

\noindent
{\it Proof}\/: Let ${\mathscr U}=\{ U_{i}\}_{i=1}^{m}$ be a covering of $Y$ by
open affine subsets and $\tilde{\mathscr U}=\{ f^{-1}(U_{i})
\}_{i=1}^{m}$. It follows from the definition of $f_{*}{\cal F}$
that $C^{0}(\tilde{\mathscr U},{\cal F})= C^{0}({\mathscr U},
f_{*}{\cal F})$. It is also obvious that $C^{p}(\tilde{\mathscr U},
{\cal F})= C^{p}({\mathscr U},f_{*}{\cal F})$ for $p>0\,$. Therefore,
$\check{H}^{p}(\tilde{\mathscr U},{\cal F})= \check{H}^{p}({\mathscr U},
f_{*}{\cal F})$. Since the \v{C}ech cohomology groups with respect to
an affine covering are isomorphic to the cohomology groups, we have
$H^{p}(X,{\cal F})\simeq H^{p}(Y,f_{*}{\cal F})$.\\

It follows from our discussion in the section on the homogeneous coordinate
ring approach that $\pi_{*}{\cal O}_{ {\por {\scriptstyle U}}_{\scriptscriptstyle\Sigma}}
=\bigoplus_{\por {\scriptscriptstyle p}}{\cal O}_{{\Bbb P}_{\scriptscriptstyle\Sigma}}({\por p})$.
Consequently, if we calculate  $H^{p}({\por U}_{\scriptscriptstyle\Sigma},{\cal O}
_{ {\por {\scriptstyle U}}_{\scriptscriptstyle\Sigma}})$ using (\ref{ck}), then 
the above lemma yields the desired result!\\


We now go on to consider  the question how an object like 
$\mbox{Ext}^{i}_{R}(M,N)$ can be algorithmically calculated.
It should be clear that the basic algorithms being 
discussed here can be equally applied to the calculations
involving Koszul complexes.  
In what follows we assume that $R={\Bbb C}[x_{1},\ldots , x_{n}]$ 
and all $R$-modules are finitely generated.
Let $M$ be an $R$-module and 
\begin{eqnarray}
  0\to K\stackrel{i}{\longrightarrow} R^{p}
        \stackrel{\pi}{\longrightarrow} M\to 0\label{ext1} 
\end{eqnarray}
a presentation thereof. The $R$-module K is said to be a (first) syzygy 
module of $M$. 
The algorithmic calculation of syzygy modules, which we will sketch 
below, is the heart of the computational approach to many questions of 
commutative algebra (cf. \cite{cox921} for more details).\\
 
Let $M=\langle \por{f}_{1},\ldots,\por{f}_{p}\rangle$ 
as a submodule of $R^m$ and $\{\por{g}_{1},\ldots,\por{g}_{q}\}$ 
be a Gr\"obner basis for $M$ such that  the leading coefficient of 
each  $\por{g}_{i}$ is one. Let $\por{x}_{i}:=\mbox{lm}(\por{g}_{i})$ be
the leading monomial of $\por{g}_{i}$ and
$\por{x}_{ij}:=\mbox{lcm}(\por{x}_{i},\por{x}_{j})$ be the least common
multiple of $\por{x}_{i}$ and $\por{x}_{j}\;$. Further, let
\begin{eqnarray}
    S(\por{g}_{i},\por{g}_{j}):=\frac{\por{x}_{ij}}{\por{x}_{i}}\;\por{g}_{i}
                              -\frac{\por{x}_{ij}}{\por{x}_{j}}\;\por{g}_{j}
                               \;.\nonumber
\end{eqnarray}
It can be shown that $S(\por{g}_{i},\por{g}_{j})$ belongs to $M$ and
therefore can be written as $S(\por{g}_{i},\por{g}_{j})=\sum_{k=1}^{q}
h_{ij}^{k}\;\por{g}_{k}$. For $i\neq j$ we define 
\begin{eqnarray}
    \por{s}_{ij}:=\frac{\por{x}_{ij}}{\por{x}_{i}}\;e_{i}
                -\frac{\por{x}_{ij}}{\por{x}_{j}}\;e_{j}
                -(h_{ij}^{1},\ldots ,h_{ij}^{q})\;,\nonumber
\end{eqnarray}
where $\{e_{i}\}_{i=1}^{q}$ is the canonical basis of $R^{q}$. The two
generating systems $F=(\por{f}_{1},\ldots,\por{f}_{p})$ and 
$G=(\por{g}_{1},\ldots,\por{g}_{q})$ of $M$ are related through the
relations $F=G\cdot \frak{A}$ and $G=F\cdot\frak{B}$, where $\frak{A}$
and $\frak{B}$ are
$q\times p$ and $p\times q$ matrices over $R$, respectively. The
matrix $\frak{A}$ can be obtained using the division algorithm. While, 
the matrix $\frak{B}$ is obtained during the calculation of the
Gr\"obner basis using the Buchberger's algorithm. Let $\por{r}_{i}\;,\;
i=1,\ldots, p$~, denote the columns of the matrix 
$1\hspace{-1mm}\mbox{l} - \frak{B}\cdot\frak{A}$. Then we have \\

\noindent
{\bf Theorem F}: $\{ \por{r}_{1},\ldots, \por{r}_{p} \}\cup
\{\frak{B}\por{s}_{ij}\mid 1\leq i<j\leq q\;\}$ generates the syzygy 
module $K=\mbox{Syz}(\por{f}_{1},\ldots,\por{f}_{p})$.\\

Using the above theorem we can easily compute the kernel of a
surjective homomorphism 
\begin{eqnarray}
  \phi :R^{p}\to M/N\;,\label{map} 
\end{eqnarray}
where $M=\langle 
\por{f}_{1},\ldots,\por{f}_{p}\rangle$ and $N=\langle \por{f}_{p+1},
\ldots,\por{f}_{p+q}\rangle$ are submodules of $R^{m}$ and
$\phi( e_{i} ) =\por{f}_{i}+N$ for $i=1,\ldots, p\;$. For that purpose
we only need to calculate $\mbox{Syz}(\por{f}_{1},\ldots,\por{f}_{p+q})$.
Let $\,\mbox{Syz}(\por{f}_{1},\ldots,\por{f}_{p+q})$ $=\langle
\por{p}_{1},\ldots,\por{p}_{\ell}\rangle$ and $\por{h}_{i}$ denote 
the vector of the first $p$ components of $\por{p}_{i}\;$. Then 
$\mbox{Ker}\;\phi=\langle\por{h}_{1},\ldots,\por{h}_{\ell}\rangle\;$.\\ 

Having the necessary computational tools at our disposal, we are now 
faced with the computation of $\mbox{Hom}(M,N)$ as a first step 
in the calculation of the Ext functor (cf. appendix). We begin with a
presentation of $M$ and $N$ 
\begin{eqnarray} 
  0\to K_{0}\stackrel{\imath_{0}}{\longrightarrow} R^{p}
   \stackrel{\pi_{0}}{\longrightarrow}  M\to 0
  & \hspace{2mm},\hspace{2mm}&
  0\to L_{0}\stackrel{\imath'_{0}}{\longrightarrow} R^{k}
   \stackrel{\pi'_{0}}{\longrightarrow} N\to 0\;.\label{ext2}
\end{eqnarray}
Similarly, we can write down presentations of the first syzygy
modules $K_{0}$ and $L_{0}$
\begin{eqnarray} 
  0\to K_{1}\stackrel{\imath_{1}}{\longrightarrow} R^{q}
   \stackrel{\pi_{1}}{\longrightarrow} K_{0}\to 0
  & \hspace{2mm},\hspace{2mm}&
  0\to L_{1}\stackrel{\imath'_{1}}{\longrightarrow} R^{\ell}
   \stackrel{\pi'_{1}}{\longrightarrow} L_{0}\to 0\;.\label{ext3}
\end{eqnarray}
Combining (\ref{ext2}) and  (\ref{ext3}) we obtain
\begin{eqnarray} 
  R^{q}\stackrel{\alpha}{\longrightarrow} R^{p}
   \stackrel{\pi_{0}}{\longrightarrow}  M\to 0
  & \hspace{2mm},\hspace{2mm}&
  R^{\ell}\stackrel{\beta}{\longrightarrow} R^{k}
   \stackrel{\pi'_{0}}{\longrightarrow} N\to 0\;,\label{ext4}
\end{eqnarray}
where $\alpha = \imath_{0}\circ\pi_{1}$ and $\beta = \imath'_{0}
\circ\pi'_{1}\,$. By considering $K_{0}$ (resp. $L_{0}$) as a submodule 
of $R^{p}$ (resp. $R^{k}$) $\alpha$ (resp. $\beta$) can be
identified with  a $p\times q$ (resp. $k\times \ell\,$) matrix whose
columns are the generators of $K_{0}$ (resp. $L_{0}$). (Note that 
theorem F provides us with such a generating system.) Applying
the (left exact contravariant functor) $\mbox{Hom}(\;\cdot\;, N)$ to
the first sequence in (\ref{ext4}) we get
\begin{eqnarray}
  0\to\mbox{Hom}(M,N)\to\mbox{Hom}(R^{p},N)
   \stackrel{\alpha^{*}}{\longrightarrow}\mbox{Hom}(R^{q},N)\;.\label{ext5}
\end{eqnarray}
Therefore, the calculation of $\mbox{Hom}(M,N)$ has been reduced to
that of the kernel of the map $\alpha^{*}:\mbox{Hom}(R^{p},N)
\to \mbox{Hom}(R^{q},N)\,,\, \psi\mapsto \alpha^{*}(\psi):=
\psi\circ \alpha\,$. For the maps (\ref{map}) we know how it can be
calculated. Hence, we are led to find  presentations of modules
$\mbox{Hom}(R^{p},N)$ and $\mbox{Hom}(R^{q},N)$. For that purpose,
we apply the exact functors $\mbox{Hom}(R^{p},\;\cdot\;)$ and
$\mbox{Hom}(R^{q},\;\cdot\;)$ to the second sequence in (\ref{ext4})
to get
\begin{eqnarray}
  \mbox{Hom}(R^{p},R^{\ell})\stackrel{\beta_{*}}{\longrightarrow}
  \mbox{Hom}(R^{p},R^{k})\to\mbox{Hom}(R^{p},N)\to 0\;,\label{ext6}\\
  \mbox{Hom}(R^{q},R^{\ell})\stackrel{\tilde{\beta}_{*}}{\longrightarrow}
  \mbox{Hom}(R^{q},R^{k})\to\mbox{Hom}(R^{q},N)\to 0\;.\label{ext7}
\end{eqnarray}
From this we see that the cokernels of the maps $\beta_{*}$ and 
$\tilde{\beta}_{*}$ yield the desired presentations. Now by calculating
the kernel of the map
\begin{eqnarray} 
    \tilde{\alpha}^{*}:\mbox{Hom}(R^{p},R^{k})\longrightarrow
    \mbox{Hom}(R^{q},R^{k})/\mbox{Im}\,\tilde{\beta}_{*}\;,\nonumber
\end{eqnarray}
where $\tilde{\alpha}^{*}$ is the composition of $\alpha^{*}$ with
the canonical projection $\mbox{Hom}(R^{p},R^{k})\to
\mbox{Hom}(R^{q},R^{k})/$
$\mbox{Im}\,\beta_{*}\,$ and taking its quotient by $\mbox{Im}\,\beta_{*}$
we arrive at $\mbox{Hom}(M,N)\,$, i.e. 
\begin{eqnarray*} 
    \mbox{Hom}(M,N)=\mbox{Ker}\,\tilde{\alpha}^{*}/\mbox{Im}\,\beta_{*}\;.
\end{eqnarray*}
We now turn to the calculation of $\mbox{Ext}^{i}_{R}(M,N)$. According
to the definition of  $\mbox{Ext}^{i}_{R}(M,N)$ we first
compute a free resolution of $M\,$, which amounts to the repeated application 
of the `Theorem F' for the calculation of syzygy modules. Note that
this procedure terminates in at most $n$ steps due to Hilbert's 
syzygy theorem.
\begin{eqnarray}
\mbox{
\unitlength0.25cm
\begin{picture}(25,14)
   \put(-9.3,10.9){\mbox{$\ldots$}}
   \put(-7,11){\vector(1,0){3.5}}
   \put(-3,10.9){\mbox{$F_{2}$}}   
   \put(2,11.7){\mbox{$\scriptscriptstyle d_{2}$}}
   \put(-.8,11){\vector(1,0){6.5}}
   \put(6,10.9){\mbox{$F_{1}$}}   
   \put(11,11.7){\mbox{$\scriptscriptstyle d_{1}$}}
   \put(8.2,11){\vector(1,0){6.5}}
   \put(15.1,10.9){\mbox{$F_{0}$}}
   \put(18.2,11.7){\mbox{$\scriptscriptstyle d_{0}$}}
   \put(17.3,11){\vector(1,0){3.5}}
   \put(21.3,10.9){\mbox{$M$}}
   \put(23.5,11){\vector(1,0){2.5}}
   \put(26.6,10.6){\mbox{$0$}}
   \put(-1.6,10){\vector(1,-1){3.5}}
   \put(.6,8.3){\mbox{$\scriptscriptstyle \pi_{2}$}}
   \put(2.8,6.5){\vector(1,1){3.5}}
   \put(5.2,8.2){\mbox{$\scriptscriptstyle \imath_{1}$}}
   \put(7.4,10){\vector(1,-1){3.5}}
   \put(9.8,8.2){\mbox{$\scriptscriptstyle \pi_{1}$}}
   \put(11.9,6.5){\vector(1,1){3.5}}
   \put(14.5,8.2){\mbox{$\scriptscriptstyle \imath_{0}$}}
   \put(1.6,4.7){\mbox{$K_{1}$}}
   \put(10.7,4.7){\mbox{$K_{0}$}}
   \put(-1.6,.5){\vector(1,1){3.5}}
   \put(3,4){\vector(1,-1){3.5}}
   \put(7.5,.5){\vector(1,1){3.5}}
   \put(12.1,4){\vector(1,-1){3.5}}
   \put(-2.5,-1){\mbox{$0$}}   
   \put(6.6,-1){\mbox{$0$}}   
   \put(15.8,-1){\mbox{$0$}}
\end{picture}}\label{free}
\end{eqnarray}

\vspace{5mm}

Next we apply the functor $\mbox{Hom}(\;\cdot\;, N)$ to (\ref{free}) to obtain
\begin{eqnarray}
  0\to\mbox{Hom}(M,N)\stackrel{d_{0}^{*}}{\longrightarrow}\mbox{Hom}(F_{0},N)
   \stackrel{d_{1}^{*}}{\longrightarrow}\mbox{Hom}(F_{1},N)
   \to \ldots \;.\label{ext8}
\end{eqnarray}
Following the steps discussed above for  the calculation of the Hom
functor we can find a presentation for each term of the sequence 
(\ref{ext8}) and then, by calculating the kernel of $d^{*}_{i}$ 
(cf. (\ref{map})) and the image of $d^{*}_{i-1}\,$, 
we get the desired result!\\



\section{Examples}\label{examples}
\hrule

\vspace{10mm} 

In the following examples we will calculate the cohomology groups
$H^{\bullet}({\por U}_{\scriptscriptstyle\Sigma},{\cal O}_
{{\por {\scriptstyle U}}_{\scriptscriptstyle\Sigma}})$ from which
$H^{\bullet}({\Bbb P}_{\scriptscriptstyle\Sigma}, {\cal O}_{{\Bbb P}_
{\scriptscriptstyle\Sigma}}({\por p}))$ can be obtained through
the projection map $\pi$ as described above.\\ 

\noindent
{\bf Example 1} (the weighted projective space $\,{\Bbb P}(w_{0},
\ldots,w_{d})\,$) \cite{dol81}:
Let $w_{0},\ldots,w_{d}$ be relatively prime
positive integers. Further, let $\{\,{\por v}_{0},\ldots,
{\por v}_{d}\,\}$ be a spanning set of a $d$-dimensional 
real vector space $V$ satisfying the linear relation 
$\;w_{0}{\por v}_{0}+\ldots+w_{d}{\por v}_{d}=0\,$. 
Let the integer span of ${\por v}_{0},\ldots,
{\por v}_{d}$ define the lattice ${\bf N}$ whose ${\Bbb R}$-extension
is obviously $V$. The fan $\Sigma$ consists of all simplicial 
cones generated by  proper subsets of $\;\{\, {\por v}_{0},\ldots,
{\por v}_{d}\,\}\,$. The corresponding toric variety will be
the weighted projective space ${\Bbb P}(w_{0},\ldots,w_{d})\,$.
It follows from these data that $S={\Bbb C}[x_{0},\ldots,x_{d}]$
with $\mbox{\sl deg}\,x_{i}=w_{i}$ and 
$I=\langle x_{0},\ldots,x_{d}\rangle\,$. 
Since the sequence  $(x_{0},\ldots,x_{d})$ is 
regular in $S={\Bbb C}[x_{0},\ldots,x_{d}]$ we find according
to theorem D and (\ref{ck}) that the only nontrivial cohomology 
groups are $H^{0}({\por U}_{\scriptscriptstyle\Sigma},
{\cal O}_{ {\por {\scriptstyle U}}_{\scriptscriptstyle\Sigma}})$ 
and $H^{d}({\por U}_{\scriptscriptstyle\Sigma},{\cal O}_
{ {\por {\scriptstyle U}}_{\scriptscriptstyle\Sigma}})\,$:
 
\begin{eqnarray}
H^{p}({\por U}_{\scriptscriptstyle\Sigma},{\cal O}_{ {\por {\scriptstyle U}}_{\scriptscriptstyle\Sigma}})=
\left\{
\begin{array}{lcc}
\hspace{4mm}S                                      &\mbox{for} & p=0    \\
                                       &           &        \\
 \hspace{4mm}                 0        &\mbox{for} & 0<p<d\\
                                       &           &        \\
     \hspace{-.5mm}\bigoplus \limits_{{n_{i}>0}\atop {i=0,\ldots,d}}
     \hspace{-2mm}{\Bbb C}\cdot x^{-n_0}_0\hspace{-2mm}
     \ldots\; x^{-n_{d}}_{d}&\mbox{for} & p=d \\
\end{array}
\right. \;  \nonumber
\end{eqnarray}
where ${\por U}_{\scriptscriptstyle\Sigma}={\Bbb A}^{d+1}\setminus \{0\}=
\bigcup_{i=0}^{d}\,U_{i}\,$. We have also made use of the 
fact that the inductive limit commutes with the cohomology (cf. 
the next example for the details on the calculation of inductive limit).\\


\noindent
{\bf Example 2}: Let ${\Bbb P}_{\scriptscriptstyle\Sigma}$ denote 
the Fano toric variety defined by the reflexive polytope $\Delta$ in $\bf N$   
with the vertices
\begin{eqnarray*}
\begin{array}{lll}
   e_{1}=(1,0,0,0,0) & e_{2}=(0,1,0,0,0) & e_{3}=(0,0,1,0,0)\hspace{4mm}
   e_{4}=(0,0,0,1,0)\\
   e_{5}=(0,0,0,0,1) & e_{6}=(0,0,0,0,-1) &
   e_{7}=(-1,-1,-3,-3,-6)\,,
\end{array}
\end{eqnarray*}
which is a  blowup of the weighted projective space ${\Bbb P}(1,1,1,3,3,6)$.
(The one-dimensional cone $\langle e_{6}\rangle$ corresponds to
the resulting exceptional divisor in ${\Bbb P}_{\scriptscriptstyle\Sigma}\,$.)
The big cones of its defining  simplicial fan $\Sigma$ are given by
\begin{eqnarray*}
\begin{array}{lllll}
&   \sigma_{1}=\langle e_{1}e_{2}e_{3}e_{4}e_{5}\rangle 
&   \sigma_{2}=\langle e_{1}e_{2}e_{3}e_{4}e_{6}\rangle
&   \sigma_{3}=\langle e_{1}e_{2}e_{3}e_{6}e_{7}\rangle
&   \sigma_{4}=\langle e_{1}e_{2}e_{4}e_{6}e_{7}\rangle\,\\
&   \sigma_{5}=\langle e_{1}e_{3}e_{4}e_{6}e_{7}\rangle
&   \sigma_{6}=\langle e_{2}e_{3}e_{4}e_{6}e_{7}\rangle
&   \sigma_{7}=\langle e_{1}e_{2}e_{3}e_{5}e_{7}\rangle
&   \sigma_{8}=\langle e_{1}e_{2}e_{4}e_{5}e_{7}\rangle\,\,\\
&   \sigma_{9}=\langle e_{1}e_{3}e_{4}e_{5}e_{7}\rangle
&   \sigma_{10}=\langle e_{2}e_{3}e_{4}e_{5}e_{7}\rangle\,.&&
\end{array}
\end{eqnarray*}
\vspace{2mm}

It follows from the above data that $S={\Bbb C}[x_{1},\ldots,x_{7}]$
with $\mbox{\sl deg}\,x_{1}=(1,0)$, $\mbox{\sl deg}\,x_{2}=(1,0)$, $\mbox{\sl deg}\,x_{3}=(3,0)$,
$\mbox{\sl deg}\,x_{4}=(3,0)$, $\mbox{\sl deg}\,x_{5}=(6,1)$, $\mbox{\sl deg}\,x_{6}=(0,1)$, 
$\mbox{\sl deg}\,x_{7}=(1,0)$, and 
$I=\langle\; x_6 x_7\,,\,x_4x_6\,,\,x_3 x_6\,,\,x_2x_6\,$, $\,
x_1x_6\,,\,x_5x_7\,,\,x_4x_5\,,\,x_3x_5\,,\,
x_2x_5\,,\,x_1x_5\;\rangle\;$.
Using the algorithm described above we get a free resolution of 
$S/I^m$ :
  $ \;0\to S\stackrel{d_6}{\longrightarrow}S^{7}\stackrel{d_5}{\longrightarrow}
    S^{20}\stackrel{d_4}{\longrightarrow}S^{30}\stackrel{d_3}{\longrightarrow}
    S^{25}\stackrel{d_2}{\longrightarrow}S^{10}\stackrel{d_1}{\longrightarrow}
    S\to S/I^m\to 0\,$ with  \\

{\scriptsize
\hspace{14mm}
$d_1=\left(
\begin{array}{c}
\mbox{\tt\hspace{-5mm}
\begin{tabular*}{10.8cm}{c@{\extracolsep\fill}c@{\extracolsep\fill}
                       c@{\extracolsep\fill}c@{\extracolsep\fill}
                       c@{\extracolsep\fill}c@{\extracolsep\fill}
                       c@{\extracolsep\fill}c@{\extracolsep\fill}
                       c@{\extracolsep\fill}c@{\extracolsep\fill}}
$x^m_6$$x^m_7$&$x^m_4$$x^m_6$&$x^m_3$$x^m_6$&$x^m_2$$x^m_6$&
$x^m_1$$x^m_6$&$x^m_5$$x^m_7$&$x^m_4$$x^m_5$&$x^m_3$$x^m_5$&
$x^m_2$$x^m_5$&$x^m_1$$x^m_5$
\end{tabular*}}
\end{array}
\hspace{-2mm}\right)$}

\vspace{6mm}

{\scriptsize
\hspace*{-4.5mm}$d_2=\left(
\begin{array}{c}
\mbox{\tt\hspace{-4.5mm}
\begin{tabular*}{14cm}{c@{\extracolsep\fill}c@{\extracolsep\fill}
                         c@{\extracolsep\fill}c@{\extracolsep\fill}
                         c@{\extracolsep\fill}c@{\extracolsep\fill}
                         c@{\extracolsep\fill}c@{\extracolsep\fill}
                         c@{\extracolsep\fill}c@{\extracolsep\fill}
                         c@{\extracolsep\fill}c@{\extracolsep\fill}
                         c@{\extracolsep\fill}c@{\extracolsep\fill}
                         c@{\extracolsep\fill}c@{\extracolsep\fill}
                         c@{\extracolsep\fill}c@{\extracolsep\fill}
                         c@{\extracolsep\fill}c@{\extracolsep\fill}
                         c@{\extracolsep\fill}c@{\extracolsep\fill}
                         c@{\extracolsep\fill}c@{\extracolsep\fill}
                         c@{\extracolsep\fill}}
 $x^m_4$&-$x^m_3$& 0&-$x^m_2$& 0& 0&-$x^m_1$& 0& 0& 0&-$x^m_5$& 0& 0&
 0& 0& 0& 0& 0&0& 0& 0& 0& 0& 0& 0 \\  
 $x^m_7$& 0&-$x^m_3$& 0&-$x^m_2$& 0& 0&-$x^m_1$& 0& 0& 0&-$x^m_5$& 0& 
 0& 0& 0& 0& 0& 0& 0& 0& 0& 0& 0& 0 \\ 
 0& $x^m_7$& $x^m_4$& 0& 0&-$x^m_2$& 0& 0&-$x^m_1$& 0& 0& 0&
 0&-$x^m_5$& 0& 0& 0& 0& 0& 0& 0& 0& 0& 0& 0 \\ 
 0& 0& 0& $x^m_7$& $x^m_4$& $x^m_3$& 0& 0& 0&-$x^m_1$& 0& 0& 0& 0& 0& 
 0&-$x^m_5$& 0& 0& 0& 0& 0& 0& 0& 0 \\ 
 0& 0& 0& 0& 0& 0& $x^m_7$& $x^m_4$& $x^m_3$& $x^m_2$& 0& 0& 0& 0& 0& 
 0& 0& 0& 0& 0&-$x^m_5$& 0& 0& 0& 0 \\ 
 0& 0& 0& 0& 0& 0& 0& 0& 0& 0& $x^m_6$& 0&-$x^m_4$& 0&-$x^m_3$& 0& 0&
 -$x^m_2$& 0& 0& 0&-$x^m_1$& 0& 0& 0 \\ 
 0& 0& 0& 0& 0& 0& 0& 0& 0& 0& 0& $x^m_6$& $x^m_7$& 0& 0&-$x^m_3$& 0& 
 0&-$x^m_2$& 0& 0& 0&-$x^m_1$& 0& 0 \\ 
 0& 0& 0& 0& 0& 0& 0& 0& 0& 0& 0& 0& 0& $x^m_6$& $x^m_7$& $x^m_4$& 0& 
 0& 0&-$x^m_2$& 0& 0& 0&-$x^m_1$& 0 \\ 
 0& 0& 0& 0& 0& 0& 0& 0& 0& 0& 0& 0& 0& 0& 0& 0& $x^m_6$& $x^m_7$&
 $x^m_4$& $x^m_3$& 0& 0& 0& 0&-$x^m_1$ \\
 0& 0& 0& 0& 0& 0& 0& 0& 0& 0& 0& 0& 0& 0& 0& 0& 0& 0& 0& 0& 
 $x^m_6$& $x^m_7$& $x^m_4$& $x^m_3$& $x^m_2$ \\ 
\end{tabular*}}
\end{array}
\hspace{-1mm}\right)$}\\

\vspace{6mm}

{\scriptsize
\hspace{-20mm}
$d_3=\left(
\begin{array}{c}
\mbox{\tt\hspace{-5mm}
\begin{tabular*}{17cm}{c@{\extracolsep\fill}c@{\extracolsep\fill}
                          c@{\extracolsep\fill}c@{\extracolsep\fill}
                          c@{\extracolsep\fill}c@{\extracolsep\fill}
                          c@{\extracolsep\fill}c@{\extracolsep\fill}
                          c@{\extracolsep\fill}c@{\extracolsep\fill}
                          c@{\extracolsep\fill}c@{\extracolsep\fill}
                          c@{\extracolsep\fill}c@{\extracolsep\fill}
                          c@{\extracolsep\fill}c@{\extracolsep\fill}
                          c@{\extracolsep\fill}c@{\extracolsep\fill}
                          c@{\extracolsep\fill}c@{\extracolsep\fill}
                          c@{\extracolsep\fill}c@{\extracolsep\fill}
                          c@{\extracolsep\fill}c@{\extracolsep\fill}
                          c@{\extracolsep\fill}c@{\extracolsep\fill}
                          c@{\extracolsep\fill}c@{\extracolsep\fill}
                          c@{\extracolsep\fill}c@{\extracolsep\fill}}
 $x^m_3$& $x^m_2$& 0& 0& $x^m_1$& 0& 0& 0& 0& 0&-$x^m_5$& 0& 0& 0& 0&
 0& 0& 0& 0& 0& 0& 0& 0& 0& 0& 0& 0& 0& 0& 0 \\ 
-$x^m_4$& 0& $x^m_2$& 0& 0& $x^m_1$& 0& 0& 0& 0& 0&-$x^m_5$& 0& 0& 0& 
 0& 0& 0& 0& 0& 0& 0& 0& 0& 0& 0& 0& 0& 0& 0 \\ 
 $x^m_7$& 0& 0& $x^m_2$& 0& 0& $x^m_1$& 0& 0& 0& 0& 0&-$x^m_5$& 0& 0& 
 0& 0& 0& 0& 0& 0& 0& 0& 0& 0& 0& 0& 0& 0& 0 \\ 
 0&-$x^m_4$&-$x^m_3$& 0& 0& 0& 0& $x^m_1$& 0& 0& 0& 0& 0& 0&-$x^m_5$& 
 0& 0& 0& 0& 0& 0& 0& 0& 0& 0& 0& 0& 0& 0& 0 \\ 
 0& $x^m_7$& 0&-$x^m_3$& 0& 0& 0& 0& $x^m_1$& 0& 0& 0& 0& 0&
 0&-$x^m_5$& 0& 0& 0& 0& 0& 0& 0& 0& 0& 0& 0& 0& 0& 0 \\ 
 0& 0& $x^m_7$& $x^m_4$& 0& 0& 0& 0& 0& $x^m_1$& 0& 0& 0& 0& 0& 0& 0&
 -$x^m_5$& 0& 0& 0& 0& 0& 0& 0& 0& 0& 0& 0& 0 \\ 
 0& 0& 0& 0&-$x^m_4$&-$x^m_3$& 0&-$x^m_2$& 0& 0& 0& 0& 0& 0& 0& 0& 0& 
 0& 0& 0&-$x^m_5$& 0& 0& 0& 0& 0& 0& 0& 0& 0 \\ 
 0& 0& 0& 0& $x^m_7$& 0&-$x^m_3$& 0&-$x^m_2$& 0& 0& 0& 0& 0& 0& 0& 0& 
 0& 0& 0& 0&-$x^m_5$& 0& 0& 0& 0& 0& 0& 0& 0 \\ 
 0& 0& 0& 0& 0& $x^m_7$& $x^m_4$& 0& 0&-$x^m_2$& 0& 0& 0& 0& 0& 0& 0& 
 0& 0& 0& 0& 0& 0&-$x^m_5$& 0& 0& 0& 0& 0& 0 \\ 
 0& 0& 0& 0& 0& 0& 0& $x^m_7$& $x^m_4$& $x^m_3$& 0& 0& 0& 0& 0& 0& 0& 
 0& 0& 0& 0& 0& 0& 0& 0& 0&-$x^m_5$& 0& 0& 0 \\
 0& 0& 0& 0& 0& 0& 0& 0& 0& 0& $x^m_4$& $x^m_3$& 0& 0& $x^m_2$& 0& 0& 
 0& 0& 0& $x^m_1$& 0& 0& 0& 0& 0& 0& 0& 0& 0 \\ 
 0& 0& 0& 0& 0& 0& 0& 0& 0& 0&-$x^m_7$& 0& $x^m_3$& 0& 0& $x^m_2$& 0& 
 0& 0& 0& 0& $x^m_1$& 0& 0& 0& 0& 0& 0& 0& 0 \\ 
 0& 0& 0& 0& 0& 0& 0& 0& 0& 0& $x^m_6$& 0& 0& $x^m_3$& 0& 0& $x^m_2$& 
 0& 0& 0& 0& 0& $x^m_1$& 0& 0& 0& 0& 0& 0& 0 \\ 
 0& 0& 0& 0& 0& 0& 0& 0& 0& 0& 0&-$x^m_7$&-$x^m_4$& 0& 0& 0& 0&
 $x^m_2$& 0& 0& 0& 0& 0& $x^m_1$& 0& 0& 0& 0& 0& 0 \\ 
 0& 0& 0& 0& 0& 0& 0& 0& 0& 0& 0& $x^m_6$& 0&-$x^m_4$& 0& 0& 0& 0& 
 $x^m_2$& 0& 0& 0& 0& 0& $x^m_1$& 0& 0& 0& 0& 0 \\ 
 0& 0& 0& 0& 0& 0& 0& 0& 0& 0& 0& 0& $x^m_6$& $x^m_7$& 0& 0& 0& 0& 0& 
 $x^m_2$& 0& 0& 0& 0& 0& $x^m_1$& 0& 0& 0& 0 \\ 
 0& 0& 0& 0& 0& 0& 0& 0& 0& 0& 0& 0& 0& 0&-$x^m_7$&-$x^m_4$& 0&
 -$x^m_3$& 0& 0& 0& 0& 0& 0& 0& 0& $x^m_1$& 0& 0& 0 \\ 
 0& 0& 0& 0& 0& 0& 0& 0& 0& 0& 0& 0& 0& 0& $x^m_6$& 0&-$x^m_4$& 0&
 -$x^m_3$& 0& 0& 0& 0& 0& 0& 0& 0& $x^m_1$& 0& 0 \\ 
 0& 0& 0& 0& 0& 0& 0& 0& 0& 0& 0& 0& 0& 0& 0& $x^m_6$& $x^m_7$& 0& 0&
 -$x^m_3$& 0& 0& 0& 0& 0& 0& 0& 0& $x^m_1$& 0 \\ 
 0& 0& 0& 0& 0& 0& 0& 0& 0& 0& 0& 0& 0& 0& 0& 0& 0& $x^m_6$& $x^m_7$& 
 $x^m_4$& 0& 0& 0& 0& 0& 0& 0& 0& 0& $x^m_1$ \\ 
 0& 0& 0& 0& 0& 0& 0& 0& 0& 0& 0& 0& 0& 0& 0& 0& 0& 0& 0& 0&-$x^m_7$&
 -$x^m_4$& 0&-$x^m_3$& 0& 0&-$x^m_2$& 0& 0& 0 \\ 
 0& 0& 0& 0& 0& 0& 0& 0& 0& 0& 0& 0& 0& 0& 0& 0& 0& 0& 0& 0& $x^m_6$&
 0&-$x^m_4$& 0&-$x^m_3$& 0& 0&-$x^m_2$& 0& 0 \\ 
 0& 0& 0& 0& 0& 0& 0& 0& 0& 0& 0& 0& 0& 0& 0& 0& 0& 0& 0& 0& 0& 
 $x^m_6$& $x^m_7$& 0& 0&-$x^m_3$& 0& 0&-$x^m_2$& 0 \\ 
 0& 0& 0& 0& 0& 0& 0& 0& 0& 0& 0& 0& 0& 0& 0& 0& 0& 0& 0& 0& 0& 0& 0& 
 $x^m_6$& $x^m_7$& $x^m_4$& 0& 0& 0&-$x^m_2$ \\
 0& 0& 0& 0& 0& 0& 0& 0& 0& 0& 0& 0& 0& 0& 0& 0& 0& 0& 0& 0& 0& 0& 0& 
 0& 0& 0& $x^m_6$& $x^m_7$& $x^m_4$& $x^m_3$ \\ 
\end{tabular*}}
\end{array}
\right)$}

\vspace{6mm}

{\scriptsize
\hspace{-21mm}
$\begin{array}{cc}
d_4\hspace{-1mm}=\hspace{-1mm}\left(
\begin{array}{c}
\mbox{\tt\hspace{-4mm}
\begin{tabular*}{11.5cm}{c@{\extracolsep\fill}c@{\extracolsep\fill}
                         c@{\extracolsep\fill}c@{\extracolsep\fill}
                         c@{\extracolsep\fill}c@{\extracolsep\fill}
                         c@{\extracolsep\fill}c@{\extracolsep\fill}
                         c@{\extracolsep\fill}c@{\extracolsep\fill}
                         c@{\extracolsep\fill}c@{\extracolsep\fill}
                         c@{\extracolsep\fill}c@{\extracolsep\fill}
                         c@{\extracolsep\fill}c@{\extracolsep\fill}
                         c@{\extracolsep\fill}c@{\extracolsep\fill}
                         c@{\extracolsep\fill}c@{\extracolsep\fill}}
 -$x^m_2$&-$x^m_1$&  0&  0&  0&-$x^m_5$&  0&  0&  0&  0&  0&  0&  0&
 0&  0&  0&  0&  0&  0&  0 \\  
  $x^m_3$&  0&-$x^m_1$&  0&  0&  0&-$x^m_5$&  0&  0&  0&  0&  0&  0&  
  0&  0&  0&  0&  0&  0&  0 \\  
 -$x^m_4$&  0&  0&-$x^m_1$&  0&  0&  0&-$x^m_5$&  0&  0&  0&  0&  0&  
  0&  0&  0&  0&  0&  0&  0 \\  
  $x^m_7$&  0&  0&  0&-$x^m_1$&  0&  0&  0&-$x^m_5$&  0&  0&  0&  0&  
  0&  0&  0&  0&  0&  0&  0 \\ 
  0&  $x^m_3$&  $x^m_2$&  0&  0&  0&  0&  0&  0&  0&-$x^m_5$&  0&  0&
  0&  0&  0&  0&  0&  0&  0 \\ 
  0&-$x^m_4$&  0&  $x^m_2$&  0&  0&  0&  0&  0&  0&  0&-$x^m_5$&  0&  
  0&  0&  0&  0&  0&  0&  0 \\
  0&  $x^m_7$&  0&  0&  $x^m_2$&  0&  0&  0&  0&  0&  0&  0&-$x^m_5$&
  0&  0&  0&  0&  0&  0&  0 \\
  0&  0&-$x^m_4$&-$x^m_3$&  0&  0&  0&  0&  0&  0&  0&  0&  0&  0&
 -$x^m_5$&  0&  0&  0&  0&  0 \\
  0&  0&  $x^m_7$&  0&-$x^m_3$&  0&  0&  0&  0&  0&  0&  0&  0&  0&
  0&-$x^m_5$&  0&  0&  0&  0 \\ 
  0&  0&  0&  $x^m_7$&  $x^m_4$&  0&  0&  0&  0&  0&  0&  0&  0&  0&
  0&  0&  0&-$x^m_5$&  0&  0 \\ 
  0&  0&  0&  0&  0&-$x^m_3$&-$x^m_2$&  0&  $x^m_1$&  0&  0&
  0&  0&  0&  0&  0&  0&  0&  0& 0 \\ 
  0&  0&  0&  0&  0&  $x^m_4$&  0&-$x^m_2$&  0&  0&  0&-$x^m_1$&  0& 
  0&  0&  0&  0&  0&  0&  0 \\ 
  0&  0&  0&  0&  0&-$x^m_7$&  0&  0&-$x^m_2$&  0&  0&  0&-$x^m_1$& 
  0&  0&  0&  0&  0&  0&  0 \\ 
  0&  0&  0&  0&  0&  $x^m_6$&  0&  0&  0&-$x^m_2$&  0&  0&  0&
 -$x^m_1$&  0&  0&  0&  0&  0&  0 \\ 
  0&  0&  0&  0&  0&  0&  $x^m_4$&  $x^m_3$&  0&  0&  0&  0&  0&  0&
 -$x^m_1$&  0&  0&  0&  0&  0 \\
  0&  0&  0&  0&  0&  0&-$x^m_7$&  0&  $x^m_3$&  0&  0&  0&  0&  0&
  0&-$x^m_1$&  0&  0&  0&  0 \\ 
  0&  0&  0&  0&  0&  0&  $x^m_6$&  0&  0&  $x^m_3$&  0&  0&  0&  0&
  0&  0&-$x^m_1$&  0&  0&  0 \\  
  0&  0&  0&  0&  0&  0&  0&-$x^m_7$&-$x^m_4$&  0&  0&  0&  0&  0&  0&
  0&  0&-$x^m_1$&  0&  0 \\ 
  0&  0&  0&  0&  0&  0&  0&  $x^m_6$&  0&-$x^m_4$&  0&  0&  0&  0&
  0&  0&  0&  0&-$x^m_1$&  0 \\ 
  0&  0&  0&  0&  0&  0&  0&  0&  $x^m_6$&  $x^m_7$&  0&  0&  0&  0&
  0&  0&  0&  0&  0&-$x^m_1$ \\
  0&  0&  0&  0&  0&  0&  0&  0&  0&  0&  $x^m_4$&  $x^m_3$&  0&  0&
  $x^m_2$&  0&  0&  0&  0&  0 \\ 
  0&  0&  0&  0&  0&  0&  0&  0&  0&  0&-$x^m_7$&  0&  $x^m_3$&  0&
  0&  $x^m_2$&  0&  0&  0&  0 \\ 
  0&  0&  0&  0&  0&  0&  0&  0&  0&  0&  $x^m_6$&  0&  0&  $x^m_3$&
  0&  0&  $x^m_2$&  0&  0&  0 \\ 
  0&  0&  0&  0&  0&  0&  0&  0&  0&  0&  0&-$x^m_7$&-$x^m_4$&  0&  0&
  0&  0&  $x^m_2$&  0&  0 \\ 
  0&  0&  0&  0&  0&  0&  0&  0&  0&  0&  0&  $x^m_6$&  0&-$x^m_4$&
  0&  0&  0&  0&  $x^m_2$&  0 \\ 
  0&  0&  0&  0&  0&  0&  0&  0&  0&  0&  0&  0&  $x^m_6$&  $x^m_7$&
  0&  0&  0&  0&  0&$x^m_2$ \\ 
  0&  0&  0&  0&  0&  0&  0&  0&  0&  0&  0&  0&  0&  0&-$x^m_7$&
 -$x^m_4$&0 &-$x^m_3$&  0&  0 \\ 
  0&  0&  0&  0&  0&  0&  0&  0&  0&  0&  0&  0&  0&  0&  $x^m_6$&
  0&-$x^m_4$&  0&-$x^m_3$&  0 \\ 
  0&  0&  0&  0&  0&  0&  0&  0&  0&  0&  0&  0&  0&  0&  0&  $x^m_6$&
  $x^m_7$&  0&  0&-$x^m_3$ \\
  0&  0&  0&  0&  0&  0&  0&  0&  0&  0&  0&  0&  0&  0&  0&  0&  0&
  $x^m_6$&  $x^m_7$&  $x^m_4$  
\end{tabular*}}
\end{array} \right)&
\hspace{-4mm}
\begin{array}{r}
d_5\hspace{-1mm}=\hspace{-1mm}\left(
\begin{array}{c}
\mbox{\tt\hspace{-4.5mm}
\begin{tabular*}{4.5cm}{c@{\extracolsep\fill}c@{\extracolsep\fill}
                      c@{\extracolsep\fill}c@{\extracolsep\fill}
                      c@{\extracolsep\fill}c@{\extracolsep\fill}
                      c@{\extracolsep\fill}}
 $x^m_1$&-$x^m_5$&  0&  0&  0&  0&  0\\  
-$x^m_2$&  0&-$x^m_5$&  0&  0&  0&  0\\  
 $x^m_3$&  0&  0&-$x^m_5$&  0&  0&  0\\  
-$x^m_4$&  0&  0&  0&-$x^m_5$&  0&  0\\  
 $x^m_7$&  0&  0&  0&  0&-$x^m_5$&  0\\  
 0&  $x^m_2$&  $x^m_1$&  0&  0&  0&  0\\  
 0&-$x^m_3$&  0&  $x^m_1$&  0&  0&  0\\  
 0&  $x^m_4$&  0&  0&  $x^m_1$&  0&  0\\  
 0&-$x^m_7$&  0&  0&  0&  $x^m_1$&  0\\  
 0&  $x^m_6$&  0&  0&  0&  0&  $x^m_1$ \\  
 0&  0&-$x^m_3$&-$x^m_2$&  0&  0&  0\\  
 0&  0&  $x^m_4$&  0&-$x^m_2$&  0&  0\\  
 0&  0&-$x^m_7$&  0&  0&-$x^m_2$&  0\\  
 0&  0&  $x^m_6$&  0&  0&  0&-$x^m_2$\\
 0&  0&  0&  $x^m_4$&  $x^m_3$&  0&  0\\  
 0&  0&  0&-$x^m_7$&  0&  $x^m_3$&  0\\                    
 0&  0&  0&  $x^m_6$&  0&  0&  $x^m_3$\\  
 0&  0&  0&  0&-$x^m_7$&-$x^m_4$&  0\\  
 0&  0&  0&  0&  $x^m_6$&  0&-$x^m_4$\\
 0&  0&  0&  0&  0&  $x^m_6$&  $x^m_7$   
\end{tabular*}}
\end{array} \right)\\
\\
d_6\hspace{-1mm}=\hspace{-1mm}\left(
\begin{array}{c}
\mbox{\tt\hspace{-4mm}
\begin{tabular*}{4cm}{c@{\extracolsep\fill}c@{\extracolsep\fill}
                      c@{\extracolsep\fill}c@{\extracolsep\fill}
                      c@{\extracolsep\fill}c@{\extracolsep\fill}
                      c@{\extracolsep\fill}}
$ x^m_{5}$ & -$x^m_{1}$ &$x^m_{2}$ & -$x^m_{3}$ & $x^m_{4}$ &
-$x^m_{7}$ &$x^m_{6}$                         
\end{tabular*}}
\end{array}
\right)^{\scriptscriptstyle T}.
\end{array}
\end{array}$}\\

\vspace{6mm}

Applying $\mbox{Hom}(\;\cdot\;,S)$ to $ \;0\to S\stackrel{d_6}
{\longrightarrow}S^{7}\stackrel{d_5}{\longrightarrow}S^{20}
\stackrel{d_4}{\longrightarrow}S^{30}\stackrel{d_3}
{\longrightarrow}S^{25}\stackrel{d_2}{\longrightarrow}
S^{10}\stackrel{d_1}{\longrightarrow}S\to 0\,$ and calculating 
the cohomology of the resulting complex as described above yields 
\begin{eqnarray*}
&&\mbox{Ext}^{2}_{S}(S/I^{m},S)=S/\langle\,x^{m}_{5},x^{m}_{6}\,\rangle\,,\\
&&\mbox{Ext}^{3}_{S}(S/I^{m},S)=0\,,\\
&&\mbox{Ext}^{4}_{S}(S/I^{m},S)=0\,,\\
&&\mbox{Ext}^{5}_{S}(S/I^{m},S)=S/\langle\,x_{1}^{m},x_{2}^{m},x_{3}^{m},
                                           x_{4}^{m},x_{7}^{m}\,\rangle\,,\\
&&\mbox{Ext}^{6}_{S}(S/I^{m},S)=S/\langle\,x_{1}^{m},x_{2}^{m},x_{3}^{m},
                                           x_{4}^{m}, x_{5}^{m},x_{6}^{m},
                                           x_{7}^{m}\rangle \,.
\end{eqnarray*}
Now we turn to the calculation of the inductive limits.  
We begin with
$(\{\mbox{Ext}^{2}_{S}(S/I^{m},S)\}$, $\{\tilde{f}_{mn}\})$.
The maps $\tilde{f}_{mn}$ are induced by $f_{mn}:S\to S\,,\, 1 \mapsto
(x_5x_6)^{n-m}\cdot 1\,$ which are in turn induced by the canonical
maps $S/I^m\to S/I^n$. Let ${\mathscr I}_m$ denote the $S$-module
$\langle\,x^{m}_{5},x^{m}_{6}\,\rangle\,$ and $\Tilde{\Tilde{f}}_{mn}
:=f_{mn}\!\!\restriction\!{\mathscr I}_m$. Then 
it can be easily seen that $\{{\mathscr I}_m\}_m$
together with the maps $\Tilde{\Tilde{f}}_{mn}$ is an inductive system.   
Since $\underset{\longrightarrow}{\mbox{lim}}
\;\mbox{Ext}^{2}_{S}(S\,/\,I^{m} , S)= \underset{\longrightarrow}
{\mbox{lim}}\, S/\,\underset{\longrightarrow}{\mbox{lim}}\,
{\mathscr I}_m$ we only need to consider the latter two limits. The first
one is well-known: $\underset{\longrightarrow}{\mbox{lim}}\, 
S= S_{x_5x_6}$, where $S_{x_5x_6}$
is the localization of $S$ at $x_5x_6\,$ 
(cf. \cite{ser55}, {\sl EGA}, {\bf 4}, p.19). We now show that
$\underset{\longrightarrow}{\mbox{lim}}\,
{\mathscr I}_m = S_{x_5}/S\oplus S_{x_6}/S\,$. It is not difficult to
verify that $(S_{x_5}/S\oplus S_{x_6}/S\,, \{\varphi_m\})$ with
$\varphi_m(x^m_5)=0\oplus\frac{1}{x^m_6}$ and 
$\varphi_m(x^m_6)=\frac{1}{x^m_5}\oplus 0$ (modulo $S$!) 
is a target object for $(\{{\mathscr I}_m\},\{\Tilde{\Tilde{f}}_{mn}\})$.
According to the universal property of the inductive limit there
exists a unique homomorphism $\psi$ such that the
diagram 
\begin{center}
\unitlength0.25cm
\begin{picture}(22,13)
   \put(3,11.2){\mbox{${\mathscr I}_m$}}
   \put(8.8,10.7){\mbox{$\scriptstyle\varphi_{m}$}}
   \put(6,11){\vector(3,-1){6.8}}
   \put(13.2,7.9){\mbox{$S_{x_5}/S\oplus S_{x_6}/S$}}
 \put(2.2,4.5){\mbox{$\underset{\longrightarrow}{\mbox{lim}}\,{\mathscr I}_m$}}
   \put(10,5.6){\mbox{$\scriptstyle \psi$}}
   \put(7.8,5.6){\vector(2,1){4.8}}
   \put(4,10.5){\vector(0,-1){4.5}}
   \put(2,8.3){\mbox{$\scriptstyle f_{m}$}}
\end{picture}
\end{center}\vspace{-8mm}
commutes for all $m\,$. We now prove that $\psi$ is actually an
isomorphism. Let $\beta$ be an arbitrary element of 
$S_{x_5}/S\oplus S_{x_6}/S\,$. Then $\beta$ can be written as
$\frac{g}{x^m_5}\oplus \frac{h}{x^m_6}$. Therefore, it follows from
 $\beta = \varphi_m(h\, x^m_5+g\,x^m_6)=\psi\circ 
f_m(h\, x^m_5+g\,  x^m_6)=\psi(f_m(h\,  x^m_5+g\,  x^m_6))$ that
 $\psi$ is surjective. Now, let $\alpha =h\,x^m_5+g\,x^m_6$ such that
$f_m(\alpha)$ lies in the kernel of $\psi\,$, i.e. $\psi(f_m(\alpha))=0\,$. 
Since $\psi\circ f_m=\varphi_m$ then 
$\varphi_m(\alpha)=0\,$, which means  $\frac{g}{x^m_5}\oplus \frac{h}{x^m_6}
=0 \Leftrightarrow (\,x^k_5\,g=0 \;\;\mbox{and}\;\; x^k_6\,h=0\,)\,$. 
Because of the commutativity of \\
\begin{center}
\unitlength0.25cm
\begin{picture}(22,12)
   \put(4.2,11.2){\mbox{${\mathscr I}_m$}}
   \put(8.8,10.7){\mbox{$\scriptstyle f_{m}$}}
   \put(7,11){\vector(2,-1){4.5}}
   \put(12.2,7.6){\mbox{$\underset{\longrightarrow}{\mbox{lim}}\,
                         {\mathscr I}_m  $}}
   \put(3.5,4.5){\mbox{${\mathscr I}_{m+k}$}}
   \put(8.8,5.4){\mbox{$\scriptstyle f_{m+k}$}}
   \put(7,5.5){\vector(2,1){4.5}}
   \put(5,10.5){\vector(0,-1){4.5}}
   \put(.6,8.4){\mbox{$\scriptstyle\Tilde{\Tilde{f}}_{m,m+k}$}}
\end{picture}
\end{center}\vspace{-8mm}
we have $f_m(\alpha)=f_{m+k}(\Tilde{\Tilde{f}}_{m,m+k}(\alpha))\,$. Since
$\Tilde{\Tilde{f}}_{m,m+k}(h\,x^m_5+g\,x^m_6)=x^k_6\,h\;x^{m+k}_5+
x^k_5\,g\;x^{m+k}_6=0\,$
we obtain $f_m(\alpha)=0\,$, i.e. $\mbox{Ker}\,\psi=\{0\}$.
Summarizing the above discussion we find out that
\begin{eqnarray*}
  H^{1}({\por U}_{\scriptscriptstyle\Sigma},{\cal O}_
  {{\por{\scriptstyle U}}_{\scriptscriptstyle\Sigma}})
& = &\hspace{-3mm}\bigoplus _{{k,l,m,n,r\geq 0}\atop {p,q>0}}
     \hspace{-2mm} {\Bbb C}\cdot x^k_1\,x^l_2\,x^m_3\,x^n_4\,
     x^{-p}_5 x^{-q}_6x^r_7
     \;\oplus\; {\Bbb C}\,[\,x_1,x_2,x_3,x_4,x_7\,] 
\end{eqnarray*}
as a ${\Bbb C}$-vector space.
The maps $\tilde{f}_{mn}$
in the next nontrivial inductive system
$(\{\mbox{Ext}^{5}_{S}(S/I^{m}$,
$S)\},\{\tilde{f}_{mn}\})$ are induced 
by $f_{mn}:S\to S\,,\, 1 \mapsto (x_1x_2x_3x_4x_7)^{n-m}\cdot 1\,$. 
As before the family $\{ {\mathscr I}_m\}_m$, where ${\mathscr I}_m$
now denotes the $S$-module
$\langle\,x_{1}^{m},x_{2}^{m},x_{3}^{m},
x_{4}^{m},x_{7}^{m}\,\rangle\,$, together with the maps 
$\Tilde{\Tilde{f}}_{mn}:=f_{mn}\!\!\restriction\!{\mathscr I}_m$
form an inductive system. Following the same lines of argumentation
as presented above one can easily verify that 
$\,\underset{\longrightarrow}{\mbox{lim}}\,
{\mathscr I}_m = S_{x_1x_2x_3x_4}/S\oplus S_{x_2x_3x_4x_7}/S\oplus
S_{x_1x_3x_4x_7}/S\oplus S_{x_1x_2x_4x_7}/S\oplus 
S_{x_1x_2x_3x_7}/S\,$. Taking $\underset{\longrightarrow}{\mbox{lim}}\, 
S= S_{x_1x_2x_3x_4x_7}\,$ into account we therefore obtain 
\begin{eqnarray*}
  H^{4}({\por U}_{\scriptscriptstyle\Sigma},{\cal O}_
  {{\por{\scriptstyle U}}_{\scriptscriptstyle\Sigma}})
& = &\hspace{-3mm}\bigoplus _{{k,l,m,n,r> 0}\atop {p,q\geq 0}}
     \hspace{-2mm}{\Bbb C}\cdot x^{-k}_1x^{-l}_2x^{-m}_3x^{-n}_4
     x^{p}_5\;x^{q}_6\;x^{-r}_7
     \;\oplus\; {\Bbb C}[\,x_5,x_6\,] 
\end{eqnarray*}
as a ${\Bbb C}$-vector space. Similar considerations in the case of
$(\{\mbox{Ext}^{6}_{S}(S/I^{m},S)\},\{\tilde{f}_{mn}\})$, where
$\tilde{f}_{mn}$ is now induced by  
$f_{mn}:S\to S\,,\, 1 \mapsto (x_1x_2x_3x_4x_5x_6x_7)^{n-m}\cdot 1\,$,
lead us to the following result
\begin{eqnarray*}
  H^{5}({\por U}_{\scriptscriptstyle\Sigma},{\cal O}_
  {{\por{\scriptstyle U}}_{\scriptscriptstyle\Sigma}})
& = &\hspace{-3mm}\bigoplus _{k,l,m,n,p,q,r> 0}
     \hspace{-2mm}{\Bbb C}\cdot x^{-k}_1\,x^{-l}_2\,x^{-m}_3\,x^{-n}_4\,
     x^{-p}_5\;x^{-q}_6\;x^{-r}_7
     \;. 
\end{eqnarray*}
 
\vspace{6mm}


\noindent
{\bf Example 3}: Let $\Delta$ be a reflexive polytope in $\bf N$ 
whose vertices are given by 
\begin{eqnarray*}
\begin{array}{llll}
e_{1}=(1,0,0,0,0)&e_{3}=(0,0,1,0,0)&e_{5}=(0,0,0,0,1)&e_{7}=(0,-1,-1,-1,-3)\\
e_{2}=(0,1,0,0,0)&e_{4}=(0,0,0,1,0)&e_{6}=(0,0,0,0,-1)&e_{8}=(-1,-2,-2,-3,-6)
\,.
\end{array}
\end{eqnarray*}
We take a maximal triangulation of the reflexive polytope $\Delta$, 
which leads to a simplicial fan $\Sigma$ whose big cones are defined by
\begin{eqnarray*}
\begin{array}{llll}
   \sigma_{1}=\langle e_{1}e_{2}e_{3}e_{4}e_{5}\rangle &
   \sigma_{2}=\langle e_{1}e_{2}e_{3}e_{4}e_{6}\rangle &
   \sigma_{3}=\langle e_{1}e_{2}e_{3}e_{5}e_{8}\rangle &
   \sigma_{4}=\langle e_{1}e_{2}e_{3}e_{6}e_{8}\rangle \\
   \sigma_{5}=\langle e_{1}e_{2}e_{4}e_{5}e_{7}\rangle &
   \sigma_{6}=\langle e_{1}e_{2}e_{4}e_{6}e_{7}\rangle &
   \sigma_{7}=\langle e_{1}e_{2}e_{5}e_{7}e_{8}\rangle &
   \sigma_{8}=\langle e_{1}e_{2}e_{6}e_{7}e_{8}\rangle \\
   \sigma_{9}=\langle e_{1}e_{3}e_{4}e_{5}e_{7}\rangle &
   \sigma_{10}=\langle e_{1}e_{3}e_{4}e_{6}e_{7}\rangle &
   \sigma_{11}=\langle e_{1}e_{3}e_{5}e_{7}e_{8}\rangle &
   \sigma_{12}=\langle e_{1}e_{3}e_{6}e_{7}e_{8}\rangle \\
   \sigma_{13}=\langle e_{2}e_{3}e_{4}e_{5}e_{8}\rangle &
   \sigma_{14}=\langle e_{2}e_{3}e_{4}e_{6}e_{8}\rangle &
   \sigma_{15}=\langle e_{2}e_{4}e_{5}e_{7}e_{8}\rangle &
   \sigma_{16}=\langle e_{2}e_{3}e_{6}e_{7}e_{8}\rangle \\
   \sigma_{17}=\langle e_{3}e_{4}e_{5}e_{7}e_{8}\rangle &
   \sigma_{18}=\langle e_{3}e_{4}e_{6}e_{7}e_{8}\rangle\,.& &
\end{array}
\end{eqnarray*}
The Fano toric variety $\,{\Bbb P}_{\scriptscriptstyle\Sigma}\,$ 
costructed from 
$\,\Sigma\,$ is a blowup of the weighted projective space 
$\,{\Bbb P}(1,1,2,2, 3,6)\,$. (The one-dimensional cones 
$\langle e_{6}\rangle$ and $\langle e_{7}\rangle$ correspond to
the resulting exceptional divisors in ${\Bbb P}_{\scriptscriptstyle\Sigma}\,$.)
It follows from the above data that $S={\Bbb C}[x_{1},\ldots,x_{8}]$
with $\mbox{\sl deg}\,x_{1}=(1,0,0)$, $\mbox{\sl deg}\,x_{2}=(2,1,0)$, $\mbox{\sl deg}\,x_{3}=(2,1,0)$,
$\mbox{\sl deg}\,x_{4}=(3,1,0)$, $\mbox{\sl deg}\,x_{5}=(6,3,1)$, $\mbox{\sl deg}\,x_{6}=(0,0,1)$, 
$\mbox{\sl deg}\,x_{7}=(0,1,0)$, $\mbox{\sl deg}\,x_{8}=(1,0,0)$ 
and $I=\langle\; x_{6}x_{7}x_{8}\,$, 
$\,x_{5}x_{7}x_{8}\,$, $\,x_{4}x_{6}x_{7}\,$, $\,x_{4}x_{5}x_{7}\,$, 
$\,x_{3}x_{6}x_{8}\,$, $\,x_{3}x_{5}x_{8}\,$, $\,x_{3}x_{4}x_{6}\,$, 
$\,x_{3}x_{4}x_{5}\,$, $\,x_{2}x_{6}x_{8}\,$, $\,x_{2}x_{5}x_{8}\,$, 
$x_{2}x_{4}x_{6}$, $x_{2}x_{4}x_{5}$, $x_{1}x_{6}x_{7}$,
$\,x_{1}x_{5}x_{7}\,$, $\,x_{1}x_{3}x_{6}\,$, $\,x_{1}x_{3}x_{5}\,$,
$\,x_{1}x_{2}x_{6}\,$, $\,x_{1}x_{2}x_{5}\;\rangle\;$.
Using the algorithm described above we get a free resolution of 
$S/I^m$ :
  $ \;0\to S\stackrel{d_6}{\longrightarrow}S^{8}\stackrel{d_5}{\longrightarrow}
    S^{27}\stackrel{d_4}{\longrightarrow}S^{48}\stackrel{d_3}{\longrightarrow}
    S^{45}\stackrel{d_2}{\longrightarrow}S^{18}\stackrel{d_1}{\longrightarrow}
    S\to S/I^m\to 0\,$.  
Proceeding as in the last example we get the following results
\begin{eqnarray*}
&&\mbox{Ext}^{2}_{S}(S/I^{m},S)=S/\langle\,x^{m}_{5},x^{m}_{6}\,\rangle\,,\\
&&\mbox{Ext}^{3}_{S}(S/I^{m},S)=S/\langle\,x^{m}_{2},x^{m}_{3},
  x^{m}_{7}\,\rangle\oplus S/\langle\,x^{m}_{1},x^{m}_{4},x^{m}_{8}\,
  \rangle\,,\\
&&\mbox{Ext}^{4}_{S}(S/I^{m},S)=S/\langle\,x^{m}_{2},x^{m}_{3},
  x^{m}_{5},x^{m}_{6},x^{m}_{7}\,\rangle\oplus S/\langle\,x^{m}_{1},x^{m}_{4},
  x^{m}_{5},x^{m}_{6},x^{m}_{8}\,\rangle\,,\\
&&\mbox{Ext}^{5}_{S}(S/I^{m},S)=S/\langle\,x_{1}^{m},x_{2}^{m},x_{3}^{m},
  x_{4}^{m},x_{7}^{m},x_{8}^{m}\,\rangle\,,\\
&&\mbox{Ext}^{6}_{S}(S/I^{m},S)=S/\langle\,x_{1}^{m},x_{2}^{m},x_{3}^{m},
  x_{4}^{m}, x_{5}^{m},x_{6}^{m}, x_{7}^{m} , x_{8}^{m}\rangle \,.
\end{eqnarray*}
By taking the inductive limit and using the results of the foregoing example
we arrive at
\begin{eqnarray*}
  H^{1}({\por U}_{\scriptscriptstyle\Sigma},{\cal O}_
  {{\por{\scriptstyle U}}_{\scriptscriptstyle\Sigma}})
& = &\hspace{-4mm}\bigoplus _{{k,l,m,n,r,s\geq 0}\atop {p,q> 0}}
     \hspace{-2mm}{\Bbb C}\cdot x^{k}_1\;x^{l}_2\;x^{m}_3\;x^{n}_4\;
     x^{-p}_5x^{-q}_6x^{r}_7\;x^{s}_8\;
     \;\oplus\; {\Bbb C}[\,x_1,x_2,x_3,x_4,x_7,x_8\,] \\
H^{2}({\por U}_{\scriptscriptstyle\Sigma},{\cal O}_
  {{\por{\scriptstyle U}}_{\scriptscriptstyle\Sigma}})
& = &\hspace{-3mm}\bigoplus _{{k,n,p,q,s\geq 0}\atop {l,m,r> 0}}
     \hspace{-2mm}{\Bbb C}\cdot x^{k}_1\;x^{-l}_2x^{-m}_3x^{n}_4\;
     x^{p}_5\;x^{q}_6\;x^{-r}_7x^{s}_8\;
     \;\oplus\; {\Bbb C}[\,x_1,x_4,x_5,x_6,x_8\,] \oplus\\
&   &\hspace{-3mm}\bigoplus _{{l,m,p,q,r\geq 0}\atop {k,n,s> 0}}
     \hspace{-2mm}{\Bbb C}\cdot x^{-k}_1x^{l}_2\;x^{m}_3\;x^{-n}_4
     x^{p}_5\;x^{q}_6\;x^{r}_7\;x^{-s}_8
     \;\oplus\; {\Bbb C}[\,x_2,x_3,x_5,x_6,x_7\,]\\ 
\end{eqnarray*}
\begin{eqnarray*}
H^{3}({\por U}_{\scriptscriptstyle\Sigma},{\cal O}_
  {{\por{\scriptstyle U}}_{\scriptscriptstyle\Sigma}})
& = &\hspace{-3mm}\bigoplus _{{k,n,s\geq 0}\atop {l,m,p,q,r> 0}}
     \hspace{-2mm}{\Bbb C}\cdot x^{k}_1\;x^{-l}_2x^{-m}_3x^{n}_4\;
     x^{-p}_5x^{-q}_6x^{-r}_7x^{s}_8\;
     \;\oplus\; {\Bbb C}[\,x_1,x_4,x_8\,] \oplus\\
&   &\hspace{-3mm}\bigoplus _{{l,m,r\geq 0}\atop {k,n,p,q,s> 0}}
     \hspace{-2mm}{\Bbb C}\cdot x^{-k}_1x^{l}_2\;x^{m}_3\;x^{-n}_4
     x^{-p}_5x^{-q}_6x^{r}_7\;x^{-s}_8
     \;\oplus\; {\Bbb C}[\,x_2,x_3,x_7\,]\\
H^{4}({\por U}_{\scriptscriptstyle\Sigma},{\cal O}_
  {{\por{\scriptstyle U}}_{\scriptscriptstyle\Sigma}})
& = &\hspace{-4mm}\bigoplus _{{k,l,m,n,r,s> 0}\atop {p,q\geq 0}}
     \hspace{-2mm}{\Bbb C}\cdot x^{-k}_1x^{-l}_2x^{-m}_3x^{-n}_4
     x^{p}_5\;x^{q}_6\;x^{-r}_7x^{-s}_8
     \;\oplus\; {\Bbb C}[\,x_5,x_6\,]\\
H^{5}({\por U}_{\scriptscriptstyle\Sigma},{\cal O}_
  {{\por{\scriptstyle U}}_{\scriptscriptstyle\Sigma}})
& = &\hspace{-3mm}\bigoplus _{{k,l,m,n,p,q,}\atop {r,s> 0}}
     \hspace{-2mm}{\Bbb C}\cdot x^{-k}_1x^{-l}_2x^{-m}_3x^{-n}_4
     x^{-p}_5x^{-q}_6x^{-r}_7x^{-s}_8\;.
\end{eqnarray*}
\hrule
\vspace{16mm}

\begin{center}
{\bf Acknowledgement}

\end{center}

\vspace{2mm}

\noindent
I would like to thank M. Kreuzer and P. Michor for helpful discussions
and comments. 
I would also like to thank SINGULAR research group 
at the mathematics department of the university of Kaiserslautern.
This work has been supported by the Austrian Research
Fund (FWF) under grant Nr. P10641-PHY and \"ONB under grant Nr. 6632.\\


\vspace{8mm}

\newpage

\noindent
{\Large\bf Appendix}\\
\hrule

\vspace{10mm}

\noindent
{\bf Inductive Limit}\\

\noindent
Let $(I,\prec )$ be a partially ordered set. $(I,\prec )$ is said 
to be directed if $\forall i,j \in I\; \exists\, k\in I$ such that 
$i\prec k$ and $j\prec k$.
An inductive system of sets\footnote{Analogously one can define
the same concepts in the category of groups, rings and modules.} 
$\mathscr S$ consists of a family of sets
$\{M_{i}\}_{i\in I}$ together with a family of maps $\{f_{ij}:M_{i}
\to M_{j}\}_{i\prec j}$ such that 1) $f_{ii}=\mbox{\sl id}_
{\scriptscriptstyle {M_{i}}}$ for all
$i\in I$, 2) $f_{jk}\circ f_{ij}=f_{ik}$ for all $i\prec j\prec k\;$.\\

A target object for the inductive system $\mathscr S$ is a set $M$
together with a family of maps $\{\varphi_{i}:M_{i}\to M\}_{i\in I}$ such
that 
\begin{center}
\unitlength0.25cm
\begin{picture}(16,12)
   \put(4,11.2){\mbox{$M_{i}$}}
   \put(8.8,10.7){\mbox{$\scriptstyle \varphi_{i}$}}
   \put(7,11){\vector(2,-1){4.5}}
   \put(12.2,7.6){\mbox{$M$}}
   \put(4,4.5){\mbox{$M_{j}$}}
   \put(8.8,5.4){\mbox{$\scriptstyle\varphi_{j}$}}
   \put(7,5.5){\vector(2,1){4.5}}
   \put(5,10.5){\vector(0,-1){4.5}}
   \put(2.8,8.4){\mbox{$\scriptstyle f_{ij}$}}
\end{picture}
\end{center}\vspace{-8mm}
commutes for all $i\prec j$. `The inductive limit of $\mathscr S$'
is a target object $(\underset{\longrightarrow}{\mbox{lim}}M_{i},
\{f_{i}\}_{i\in I})$ satisfying the following `universal property':
for any target object $(M,\{\varphi_{i}\}_{i\in I})$ there is a unique
map $f:\underset{\longrightarrow}{\mbox{lim}}M_{i}\to M$ such that the
diagram 
\begin{center}
\unitlength0.25cm
\begin{picture}(16,13)
   \put(3,11.2){\mbox{$M_{i}$}}
   \put(8.8,10.7){\mbox{$\scriptstyle \varphi_{i}$}}
   \put(6,11){\vector(3,-1){6.8}}
   \put(13.2,7.9){\mbox{$M$}}
   \put(2.2,4.5){\mbox{$\underset{\longrightarrow}{\mbox{lim}}M_{i}$}}
   \put(10,5.6){\mbox{$\scriptstyle f$}}
   \put(7.8,5.6){\vector(2,1){4.8}}
   \put(4,10.5){\vector(0,-1){4.5}}
   \put(2.5,8.3){\mbox{$\scriptstyle f_{i}$}}
\end{picture}
\end{center}\vspace{-12mm}
commutes for all $i\in I$.\\

\vspace{4mm}
\noindent
{\bf The Ext Functor}\\

\noindent
We first recall some elementary
concept of homological algebra. Let $R$ be a ring. A chain
complex of $R$-modules is a sequence of $R$-modules and homomorphisms 
\begin{eqnarray*}
  (M_{\bullet},d_{\bullet})\; : \;
  \ldots \to M_{i+1}\stackrel{d_{i+1}}{\longrightarrow}M_{i}
  \stackrel{d_{i}}{\longrightarrow}M_{i-1}\to\ldots
\end{eqnarray*}
with $d_{i}\circ d_{i+1}=0$ for all $i\in {\Bbb Z}$. 
A cochain complex of
$R$-modules is a sequence of $R$-modules and homomorphisms
\begin{eqnarray*}
  (M^{\bullet},d^{\bullet})\; : \;
  \ldots \to M^{i-1}\stackrel{d^{i-1}}{\longrightarrow}M^{i}
  \stackrel{d^{i}}{\longrightarrow}M^{i+1}\to\ldots
\end{eqnarray*}
with  $d^{i}\circ d^{i-1}=0$ for all $i\in {\Bbb Z}$. 
The homology of a chain complex $ (M_{\bullet},d_{\bullet})$
is given by
$H_{i}(M_{\bullet}):=\mbox{Ker}\;d_{i}/\mbox{Im}\;d_{i+1}$.
The cohomology of a cochain complex $(M^{\bullet},d^{\bullet})$
is given by $H^{i}(M^{\bullet}):=\mbox{Ker}\;d^{i}/\mbox{Im}\;d^{i-1}$.\\

Let $ (M_{\bullet},d_{\bullet})$ and $ (N_{\bullet},d'_{\bullet})$ be
two chain complexes. A homomorphism of chain complexes
$f_{\bullet}:M_{\bullet}\to N_{\bullet}$ is a family of homomorphisms 
$f_{i}:M_{i}\to N_{i}$ such that $f_{i-1}\circ d_{i}=d'_{i}\circ
f_{i}\;$ for all $\,i\in {\Bbb Z}\,$. 
Clearly $f_{\bullet}:M_{\bullet}\to N_{\bullet}$ induces
a well-defined map $H_{i}(f_{\bullet}):H_{i}(M_{\bullet})\to 
H_{i}(N_{\bullet})$ for all $i\in {\Bbb Z}$. Two homomorphisms
$f_{\bullet},g_{\bullet}:M_{\bullet}\to N_{\bullet}$ are said to be
homotopic, written $f_{\bullet}\simeq g_{\bullet}$, if there are maps
$h_{i}:M_{i}\to N_{i+1}, i\in {\Bbb Z}$, such that $ f_{i}-g_{i}=
d'_{i+1}\circ h_{i}+h_{i-1}\circ d_{i}\;$. The chain complexes
$ (M_{\bullet},d_{\bullet})$ and $ (N_{\bullet},d'_{\bullet})$ 
are called homotopy equivalent, written $ M_{\bullet}\simeq 
N_{\bullet}\,$, if there are homomorphisms $f_{\bullet}:M_{\bullet}
\to N_{\bullet}$ and $g_{\bullet}:N_{\bullet}\to M_{\bullet}$ 
such that $f_{\bullet}\circ g_{\bullet}\simeq\mbox{id}_{\scriptscriptstyle
N_{\bullet}}$ and $g_{\bullet}\circ f_{\bullet}\simeq\mbox{id}_
{\scriptscriptstyle
M_{\bullet}}$. Homotopy equivalent chain complexes have the 
same homology.\\

An $R$-module $P$ is called projective if for any surjective
homomorphism $f:M\to N$ of $R$-modules $M$ and $N$ and any
homomorphism $g:P\to N$ there exists a homomorphism $h:P\to M$
such that $g=f\circ h$. It follows from the definition that a free
module is projective.
By a projective (resp. free) resolution of an
$R$-module $M$ we mean an exact sequence 
 $ 
 \ldots \to P_{i}\stackrel{d_{i}}{\longrightarrow}P_{i-1}
 \to\ldots \to P_{0}\stackrel{d_{0}}{\longrightarrow}M\to 0$,
where $P_{i}$ are projective (resp. free) for $i=0,1,\ldots$. It
can be shown that any two projective resolutions of $M$ are
homotopy equivalent.\\

Let $M$ and $N$ be two $R$-modules and $P_{\bullet}\to M\to 0$ a 
projective resolution of $M$. Then 
\begin{eqnarray*}
\mbox{Ext}_{R}^{i}(M,N):= 
H^{i}(\mbox{Hom}(P_{\bullet},N)).
\end{eqnarray*}





\end{document}